\documentclass[11pt,french]{amsart}
\usepackage[francais]{babel}
\usepackage[utf8]{inputenc}
\usepackage{latexsym}
\usepackage[all]{xy}
\usepackage{amsmath}
\usepackage{amssymb}
\usepackage{amsfonts}

\def\N{{\mathbb N}}

\def\tf{$\sqcup\hskip-2.3mm\sqcap$}   

{\newtheorem{theorem}{Th\'eor\`eme}[section]   
   \newtheorem{proposition}[theorem]{Proposition}
   \newtheorem{lemma}[theorem]{Lemme}  
   \newtheorem{remark}[theorem]{Remarque}

{\theoremstyle{definition}   
   }
{\theoremstyle{remark}  
   \newtheorem*{definition}{D\'efinition}    
    
   \newtheorem*{preuve}{Preuve}

   \newtheorem*{preuvele}{Preuve du lemme}  

\author{Michel Vaqui\'e$^*$} 
\address{CNRS, Universit\'e de Toulouse, Institut de Math\'ematiques de Toulouse (UMR 5219),  
118 route de Narbonne, 31062 Toulouse Cedex 9, 
France} 
\email{vaquie@math.univ-toulouse.fr} 
\title[Famille admise]{Famille admise associ\'ee \`a une valuation de $\mathbf{K(X)}$}

   \thanks{$*$ Partially supported by the grant of the Agence Nationale de la Recherche “CatAG”ANR-17-CE40-0014.}

\begin{document}  

\begin{abstract} 
Soit $K$ un corps muni d'une valuation $\nu$, et soit $L=K(x)$ une extension monog\`ene transcendante de $K$, alors toute valuation $\mu$ de $L$ qui prolonge $\nu$ est d\'etermin\'ee par sa restriction \`a l'anneau des polyn\^omes $K[x]$. 
Nous savons associer \`a cette valuation $\mu$ une famille de valuations $\mathcal{A} = (\mu _i) _{i \in I}$ 
de $K[x]$, appel\'ee famille admise associ\'ee, qui converge en un certain sens vers la valuation $\mu$. 
Bien que la d\'efinition de cette famille, ainsi que la notion de convergence, fassent intervenir de mani\`ere esentielle la structure de l'anneau de polyn\^omes, en particulier le degr\'e des polyn\^omes, nous montrons dans cette note que la famille $\mathcal{A}$ de valuations de $L$ ne d\'epend pas du g\'en\'erateur $x$ choisi.

\vskip .2cm

\noindent {\scshape Abstract}. Let $K$ be a field with a valuation $\nu$ and let $L=K(x)$ be a transcendantal extension of $K$, then any valuation $\mu$ of $L$ which extends $\nu$ is determined by its restriction to the polynomial ring $K[x]$. 
We know how to associate to this valuation $\mu$ a family of valuations $\mathcal{A} = (\mu _i) _{i \in I}$ of $K[x]$, 
called the associated admise family, which converges in a certain sense towards the valuation $\mu$. 
Although the definition of this family, as well as the notion of convergence, essentially imply the structure of the polynomial ring, in particular the degree of polynomials, we show in this note that the family $\mathcal{A}$ of valuations of $L$ do not depend on the chosen generator $x$. 
 \end{abstract}  


\vskip .2cm 

\date{Octobre 2018} 

\maketitle   

\tableofcontents

		      \section*{Introduction}

Soit $K$ un corps muni d'une valuation $\nu$, \`a valeurs dans un groupe ordonn\'e $\Gamma _{\nu}$, une famille admise de valuations $\mathcal{A} = ( \mu _i) _{i \in I}$ est une famille de valuations ou de pseudo-valuations distinctes de l'anneau des polyn\^omes $K[x]$ prolongeant la valuation $\nu$, \`a valeurs dans un groupe totalement ordonn\'e $\Gamma$, index\'ee par un ensemble totalement ordonn\'e $I$ poss\'edant un plus petit \'el\'ement $\iota$, v\'erifiant les propri\'et\'es suivantes:  
\begin{enumerate}
\item pour tout polyn\^ome $f$ dans $K[x]$ la famille $( \mu _i(f)) _{i \in I}$ est croissante dans le groupe $\Gamma$, c'est-\`a-dire $\mu _i(f) \leq \mu _j (f)$ pour $i<j$ dans $I$, de plus s'il existe $i<j$ avec $\mu _i(f) = \mu _j(f)$ alors pour tout $i' \geq i$ nous avons $\mu _i(f) = \mu _{i'} (f)$; 
\item pour tout polyn\^ome $f$ dans $K[x]$ la famille $( \mu _i(f)) _{i \in I}$ admet un plus grand \'el\'ement dans le groupe $\Gamma$, c'est-\`a-dire si l'ensemble $I$ n'a pas de plus grand \'el\'ement pour tout polyn\^ome $f$ il existe $i \in I$ tel que $\mu _i(f) = \mu _j(f)$ pour tout $j \geq i$.
\end{enumerate}

Comme l'anneau $K[x]$ est de dimension un, si $\mu _1$ est une pseudo-valuation de $K[x]$ dont le \emph{socle} $\mathfrak{p} = \{ x \in K[x] \ | \ \mu _1 (x) = + \infty \}$ est non r\'eduit \`a $(0)$, il n'existe pas de pseudo-valuation $\mu _2$ de $K[x]$ distincte de $\mu _1$ telle que $\mu _1 \leq \mu _2$.   
Par cons\'equent toutes les pseudo-valuations d'une famille admise $\mathcal{A} = ( \mu _i) _{i \in I}$ sont en fait des valuations, \`a l'exception de la \emph{derni\`ere}, c'est-\`a-dire de $\mu _{\bar \iota}$ dans le cas o\`u l'ensemble $I$ admet un plus grand \'el\'ement $\bar \iota$, qui peut \^etre une valuation ou une pseudo-valuation.  

Nous appelons \emph{limite} de la famille admise $\mathcal{A} = ( \mu _i) _{i \in I}$ la valuation ou pseudo-valuation $\mu$ de $K[x]$ d\'efinie pour tout polyn\^ome $f$ par 
$$\mu (f) \ = \ Sup \bigl ( \mu _i (f) \ ; \ i \in I \bigr ) \ .$$ 
Si l'ensemble $I$ admet un plus grand \'el\'ement $\bar \iota$, la limite de la famille est la valuation ou pseudo-valuation $\mu _{\bar \iota}$, sinon la limie est une valuation d\'efinie par $\mu (f) = \mu _i(f)$ pour $i$ assez grand dans $I$. 
Nous disons aussi que la famille $\mathcal{A}$ \emph{converge} vers la valuation $\mu$. 

De plus chaque valuation $\mu _i$ de la famille est d\'efinie de mani\`ere explicite \`a partir des valuations $\mu _j$ pour $j < i$ et d'un polyn\^ome unitaire irr\'eductible $\phi _i$ de $K[x]$ appel\'e suivant les cas \emph{polyn\^ome-cl\'e} ou \emph{polyn\^ome-cl\'e limite}. 
Nous associons ainsi \`a la famille admise $\mathcal{A}$ une famille de polyn\^omes unitaires $( \phi _i) _{i \in I}$, et cette famille v\'erifie les deux propri\'et\'es suivantes: 
\begin{enumerate}
\item pour $i<j$ dans $I$ nous avons deg $\phi _i \leq$ deg $\phi _j$ et deg$\phi _{\iota} =1$ o\`u $\iota$ est le plus petit \'el\'ement de l'ensemble $I$; 
\item pour tout polyn\^ome $f$ de $K[x]$, si deg $f <$ deg $\phi _i$ alors $\mu _i (f) = \mu (f)$. 
\end{enumerate}   

\vskip .2cm 

Soit $L$ une extension monog\`ene de $K$, alors toute valuation $\tilde\mu$ de $L$ qui prolonge $\nu$ est d\'etermin\'ee par une valuation ou une pseudo-valuation $\mu$ de l'anneau des polyn\^omes $K[x]$ qui prolonge $\nu$. Plus pr\'ecis\'ement $\mu$ est une valuation de $K[x]$, restriction de la valuation $\tilde \mu$ si $L=K(x)$ est une extension transcendante pure de $K$ et $\mu$ est une pseudo-valuation de socle $\mathfrak{p}$ si $L$ est une extension alg\'ebrique de $K$ d\'efinie par un polyn\^ome $\phi$, $L = K[x] / ( \phi )$, et l'id\'eal $\mathfrak{p}$ est l'id\'eal de $K[x]$ engendr\'e par $\phi$. 

Nous savons que nous pouvons obtenir toute valuation ou pseudo-valuation $\mu$ de l'anneau des polyn\^omes $K[x]$ qui prolonge $\nu$ gr\^ace \`a une famille admise de valuations ${\mathcal A}=\bigl ( \mu _i \bigr )_{i \in I}$ qui converge vers $\mu$, et cette famille est essentiellement unique (cf. th\'eor\`emes 2.4 et 2.5 de \cite{Va 1} et proposition 2.9. de \cite{Va 2}), et nous disons que $\mathcal{A}$ est une \emph{famille admise associ\'ee} \`a la valuation ou la pseudo-valuation $\mu$ (cf. \cite{Va 2}). 

\vskip .2cm 

La d\'efinition de la famille admise $\mathcal{A} = (\mu _I)_{i \in I}$ et sa construction explicite \`a partir de la famille de polyn\^omes $(\phi _i)_{i \in I}$ utilisent de mani\`ere essentielle les propri\'et\'es de l'anneau $K[x]$, en particulier la notion de degr\'e et l'existence de la division euclidienne. 
De m\^eme la propri\'et\'e de croissance des valeurs $\mu _i(f)$ pour $f$ dans l'anneau $K[x]$ est essentielle \`a la d\'efinition de la famille admise, alors que si $\mu _1$ et $\mu _2$ sont deux valuations distinctes sur un corps $L$, il ne peut exister de relation de comparaison de la forme $\mu _1 (u) \leq \mu _2 (u)$ valable pour tout $u$ dans $L$. 

Pourtant dans le cas o\`u l'extension $L$ est l'extension transcendante pure $K(x)$, si $\mu$ est une valuation de $L$ qui prolonge la valuation $\nu$, la famille $\mathcal{A} = (\mu _i) _{i \in I}$, consid\'er\'ee comme famille de valuations de $L$, est ind\'ependante du g\'en\'erateur $x$ choisi. 
Plus pr\'ecis\'ement nous allons montrer que si $x$ et $y$ sont deux g\'en\'erateurs de $L$ sur $K$, les familles admises de valuations de $L$ associ\'ees aux restrictions de $\mu$ aux anneaux des polyn\^ome $K[x]$ et $K[y]$ sont \'equivalentes.

      \section{Famille admise associ\'ee \`a une valuation de $K[x]$}    

Dans ce qui suit nous nous donnons une valuation $\nu$ sur un corps $K$ et toutes les valuations ou pseudo-valuations $\mu$ de l'anneau des polyn\^omes $K[x]$ que nous consid\'erons sont des prolongements de $\nu$. 
Nous nous donnons aussi un groupe totalement ordonn\'e $\tilde\Gamma$, contenant le groupe des ordres $\Gamma _{\nu}$ de la valuation $\nu$, et toutes les valuations ou pseudo-valuations $\mu$ de $K[x]$ ont leur groupe des ordres $\Gamma _{\mu}$ qui est un sous-groupe ordonn\'e de $\tilde\Gamma$.   

Pour toute valuation $\mu$ de $K[x]$ nous pouvons d\'efinir la notion de {\it polyn\^ome-cl\'e} $\phi$, et si $\phi$ est un polyn\^ome-cl\'e pour $\mu$ et si $\gamma$ est un \'el\'ement de $\tilde\Gamma$ v\'erifiant $\gamma > \mu (\phi )$, nous pouvons d\'efinir une nouvelle valuation $\mu '$ de $K[x]$, appel\'ee {\it valuation augment\'ee} associ\'ee au polyn\^ome-cl\'e $\phi$ et \`a la valeur $\gamma$ que nous notons $\mu ' = [ \mu \ ; \ \mu ' (\phi ) = \gamma ]$, de la mani\`ere suivante:  
   
pour tout polyn\^ome $f$ de $K[x]$, nous \'ecrivons le d\'eveloppement de $f$ selon les puissances de $\phi$, $f = g_m \phi ^m + \ldots + g_1 \phi + g_0$, o\`u les polyn\^omes $g_j$, $0 \leq j \leq m$ sont de degr\'e strictement inf\'erieur au degr\'e du polyn\^ome-cl\'e $\phi$, et nous avons:  
$$ \mu ' (f) \ = \ Inf \left ( \mu ( g_j ) + j \gamma \ ; \ 0 \leq j \leq m \right) \ .$$   

Nous pouvons d\'efinir aussi pour tout polyn\^ome unitaire $\phi$ de degr\'e un, $\phi = x+b$, et pour toute valeur $\gamma$ de $\tilde\Gamma$ une valuation $\mu$ de $K[x]$, que nous appelons encore {\it valuation augment\'ee} associ\'ee au polyn\^ome-cl\'e $\phi$ et \`a la valeur $\gamma$ que nous notons $\mu  = [ \nu \ ; \ \mu  (\phi ) = \gamma ]$, de la mani\`ere suivante:  

tout polyn\^ome $f$ de $K[x]$ s'\'ecrit de mani\`ere unique sous la forme $f= a_d \phi ^d + \ldots + a_1 \phi + a_0$, avec $a_j \in K$, et nous posons 
$$\mu (f) \ = \ Inf \left ( \nu (a_j) + j \gamma \ ; \ 0 \leq j \leq d \right ) \ .$$

De plus nous pouvons aussi d\'efinir la notion de {\it famille de valuations augment\'ees it\'er\'ees} comme une famille d\'enombrable $\bigl ( \mu _i \bigr ) _{i \in I}$ de valuations de $K[x]$, $I = \{ 1, \ldots ,n \}$ ou $I = \N ^*$, associ\'ee \`a une famille de polyn\^omes $\bigl (\phi _i\bigr )_{i \in I}$ et \`a une famille $\bigl (\gamma _i\bigr )_{i \in I}$ d'\'el\'ements de $\tilde\Gamma$, telle que chaque valuation $\mu _i$, $i>1$, est une valuation augment\'ee de la forme $\mu _i = [ \mu _{i-1} \ ; \ \mu _i (\phi _i) = \gamma _i ]$ et o\`u la famille des polyn\^omes-cl\'es $\bigl (\phi _i \bigr )$ v\'erifie les deux propri\'et\'es suivantes: pour tout $i>2$ nous avons deg$\phi _i \geq$ deg$\phi _{i-1}$ et les polyn\^omes $\phi _i$ et $\phi _{i-1}$ ne sont pas $\mu _{i-1}$-\'equivalents.    
Nous renvoyons aux articles \cite{McL 1}, \cite{McL 2}, et \cite {Va 1}, pour les d\'efinitions et les propri\'et\'es des polyn\^omes-cl\'es, des valuations augment\'ees et des familles de valuations augment\'ees it\'er\'ees.      

\vskip 0.2cm  

De mani\`ere analogue pour toute famille admissible continue ${\mathcal C} = \bigl ( \mu _{\alpha} \bigr ) _{\alpha \in A}$ ( voir la d\'efinition pr\'ecise plus bas) nous pouvons d\'efinir la notion de {\it polyn\^ome-cl\'e limite} et de \emph{valuation augment\'ee limite}.  
Plus pr\'ecis\'ement soit ${\mathcal C} = \bigl ( \mu _{\alpha} \bigr ) _{\alpha \in A}$ une famille admissible continue de valuations de $K[x]$, index\'ee par un ensemble totalement ordonn\'e $A$ sans plus grand \'el\'ement, associ\'ee \`a la famille de polyn\^omes-cl\'es $\bigl ( \phi _{\alpha} \bigr ) _{\alpha \in A}$ et \`a la famille de valeurs $\bigl ( \gamma _{\alpha} \bigr ) _{\alpha \in A}$. 
Par d\'efinition chaque valuation $\mu _{\alpha}$ est une valuation augment\'ee de la forme $\mu _{\alpha} = [ \mu  \ ; \ \mu _{\alpha} (\phi _{\alpha}) = \gamma _{\alpha} ]$, o\`u $\mu$ est une valuation de $K[x]$ donn\'ee, les polyn\^omes-cl\'es $\phi _{\alpha}$ sont tous de m\^eme degr\'e $d$ et les valeurs $\gamma _{\alpha}$ forment une famille croissante sans plus grand \'el\'ement dans $\tilde\Gamma$.   

Nous d\'efinissons alors l'ensemble 
$$\tilde\Phi \bigl ( \mathcal{C} \bigr ) = \tilde\Phi \left ( {\bigl (\mu _{\alpha}\bigr )}_{\alpha \in A} \right ) \ = \ \{ f \in K[x] \ | \ \mu _{\alpha} (f) < \mu _{\beta}(f) , \forall \alpha < \beta \in A \} \ ,$$  
Si cet ensemble est non vide, nous appelons $d _{ \mathcal{C}}$ le degr\'e minimal des polyn\^omes appartenant \`a $\tilde\Phi _{\mu}\bigl ( \mathcal{C} \bigr )$ et nous d\'efinissons l'ensemble  
$$\Phi \bigl ( \mathcal{C} \bigr ) = \Phi \left ( {\bigl (\mu _{\alpha} \bigr )}_{\alpha \in A} \right ) \ = \ \{ \phi \in \tilde\Phi \bigl ( \mathcal{C} \bigr ),\  \hbox{deg}\phi =  d _{ \mathcal{C}} \ \hbox{et $\phi$ unitaire} \  \} \ .$$  
Un polyn\^ome $\phi$ appartenant \`a $\Phi \bigl ( \mathcal{C} \bigr )$ est un polyn\^ome-cl\'e-limite pour la famille $\mathcal{C}$, mais nous pouvons donner une autre d\'efinition de polyn\^ome-cl\'e limite de la mani\`ere suivante. 

Nous disons qu'un polyn\^ome $f$ \emph{$A$-divise} un polyn\^ome $g$ s'il existe $\alpha \in A$ tel que pour tout $\beta \geq \alpha$ dans $A$ , $f$ $\mu _{\beta}$-divise $g$. 
Nous pouvons alors d\'efinir les notions de \emph{$A$-minimalit\'e} et \emph{$A$-irr\'eductibilit\'e} et un polyn\^ome-cl\'e-limite pour la famille ${\mathcal C} = \bigl ( \mu _{\alpha} \bigr ) _{\alpha \in A}$ est un polyn\^ome unitaire $\phi$ qui est $A$-minimal et $A$-irr\'eductible. 

Soit $\phi$ un polyn\^ome-cl\'e limite $\phi$ pour la famille $\mathcal{C}$, alors un polyn\^ome $f$ appartient \`a $\tilde\Phi \bigl ( \mathcal{C} \bigr )$ si et seulement si il est $A$-divisible par $\phi$, et pour tout polyn\^ome $f$ n'appartenant pas \`a $\tilde\Phi \bigl ( \mathcal{C} \bigr )$, en particulier pour tout polyn\^ome $f$ de degr\'e deg$f<$ deg$\phi$, nous pouvons d\'efinir $\mu _A(f)$ par $\mu _A (f) = Sup \left ( \mu _{\alpha} (f) \ ; \ \alpha \in A \right ) = \mu _{\alpha}(f)$ pour $\alpha$ suffisamment grand. 

Soient $\phi$ un polyn\^ome-cl\'e limite et $\gamma$ un \'el\'ement de $\tilde\Gamma$ v\'erifiant $\gamma > \mu _{\alpha}(\phi )$ pour tout $\alpha$ dans $A$, nous pouvons d\'efinir une nouvelle valuation $\mu '$ de $K[x]$, appel\'ee {\it valuation augment\'ee limite} pour $\mathcal{C}$ associ\'ee au polyn\^ome-cl\'e limite $\phi$ et \`a la valeur $\gamma$ que nous notons $\mu ' = \bigl [ ( \mu _{\alpha} ) _{\alpha \in A} \ ; \ \mu ' (\phi ) = \gamma \bigr ]$, de la mani\`ere suivante:  

pour tout polyn\^ome $f$ de $K[x]$, nous \'ecrivons le d\'eveloppement de $f$ selon les puissances de $\phi$, $f = g_m \phi ^m + \ldots + g_1 \phi + g_0$, o\`u les polyn\^omes $g_j$, $0 \leq j \leq m$, sont de degr\'e strictement inf\'erieur au degr\'e du polyn\^ome-cl\'e limite $\phi$, et nous posons:           
$$ \mu ' (f) \ = \ Inf \left (  \mu _A ( g_j ) + j \gamma \  ; \ 0 \leq j \leq m  \right ) \ .$$    

Nous allons finir cette partie par un r\'esultat qui est l'analogue pour les valuations augment\'ees limites du th\'eor\`eme 5.2. de \cite{McL 1}. Pour cela nous devons rappeler les propri\'et\'es des polyn\^omes-cl\'es limites et des valuations augment\'ees limites. 

\begin{lemma}
Soit $f=q \phi +r$ la division euclidienne de $f$ par le polyn\^ome-cl\'e limite $\phi$, alors si $f$ n'est pas $A$-divisible par $\phi$ nous avons $\mu _A(r) = \mu _A (f) < \mu _{\alpha}(q\phi )$ pour tout $\alpha$ suffisamment grand dans $A$. 
\end{lemma}

\begin{preuvele} 
Comme les polyn\^omes $r$ et $f$ n'appartiennent pas \`a l'ensemble $\tilde\Phi \bigl ( \mathcal{C} \bigr )$ il existe $\alpha$ dans $A$ tel que pour tout $\beta > \alpha$ nous ayons $\mu _{\beta} (r) = \mu _{\alpha}(r) = \mu _A(r)$ et $\mu _{\beta} (f) = \mu _{\alpha}(f) = \mu _A(f)$. 
Nous d\'eduisons le r\'esultat des in\'egalit\'es strictes $\mu _{\alpha} (q\phi ) < \mu _{\beta}(q\phi )$ pour tout $\alpha < \beta$ dans $A$. 
 
\hfill \tf    
\end{preuvele} 

\begin{proposition}\label{prop-definition-val-augm-lim} 
Soit $\phi$ un polyn\^ome-cl\'e limite pour la famille continue $\mathcal C$ et soit $\mu '$ la valuation augment\'ee limite associ\'ee \`a $\phi $ et \`a une valeur $\gamma$, $\mu ' = \bigl [ ( \mu _{\alpha} ) _{\alpha \in A} \ ; \ \mu ' (\phi ) = \gamma \bigr ]$. 
Pour tout polyn\^ome $f$ de $K[x]$, si nous avons une relation $f = f_m \phi ^m + \ldots + f_1 \phi + f_0$ avec des polyn\^omes $f_j$, $0 \leq j \leq m$, non $A$-divisibles par $\phi$ nous avons l'\'egalit\'e 
$$ \mu ' (f) \ = \ Inf \left (  \mu _A ( f_j ) + j \gamma \  ; \ 0 \leq j \leq m  \right ) \ .$$   
\end{proposition} 

\begin{preuve} 
Pour tout $j$ soit $f_j=q_j \phi +r_j$ la division euclidienne de $f_j$ par $\phi$, et nous posons $f=q \phi + r$ avec $r= r_m \phi ^m + \ldots + r_1 \phi + r_0$ et $q = q_m \phi ^m + \ldots + q_1 \phi + q_0$. 
Nous d\'eduisons du lemme pr\'ec\'edent 
$$\begin{array}{rcl}
\mu ' (f) & \geq & Inf \left (  \mu _A ( f_j ) + j \gamma \  ; \ 0 \leq j \leq m  \right ) \\ 
& = & Inf \left (  \mu _A ( r_j ) + j \gamma \  ; \ 0 \leq j \leq m  \right )  \ = \ \mu '(r) \ ,
\end{array}$$  
et pour $\alpha$ suffisamment 
$$\begin{array}{rcl} 
\mu ' (q \phi) & \geq & Inf \left (  \mu '( q_j \phi) + j \gamma \  ; \ 0 \leq j \leq m  \right ) \\  
& > & Inf \left (  \mu _{\alpha} ( q_j \phi) + j \gamma \  ; \ 0 \leq j \leq m  \right ) \\
& > & Inf \left (  \mu _A ( r_j ) + j \gamma \  ; \ 0 \leq j \leq m  \right ) \ = \ \mu '(r) \ , 
\end{array}$$
d'o\`u l'\'egalit\'e cherch\'ee.  
 
\hfill \tf    
\end{preuve} 

Nous renvoyons \`a \cite{Va 1} pour les d\'efinitions pr\'ecises et les propri\'et\'es des poly\-n\^omes-cl\'es limites et des valuations augment\'ees limites.

\vskip .2cm 

Nous pouvons maintenant d\'efinir une famille admissible, et aussi fixer les notations que nous utiliserons dans la suite.  

\begin{definition}  
Une {\it famille admissible $\mathcal A$ pour la valuation $\nu$ de $K$} est une famille de valuations $\bigl ( \mu _i \bigr ) _{i \in I}$ de $K[x]$, obtenue comme r\'eunion de familles admissibles simples    
$${\mathcal A} \ = \ \bigcup _{j \in J} {\mathcal S}^{(j)} \ , $$  
o\`u $J$ est un ensemble d\'enombrable, $J=\{ 1, \ldots , N\}$ ou $J = \N^*$, et nous d\'efinissons $J^*$ par $J^*=\{ 1, \ldots , N-1\}$ si $J$ est fini et par $J^*=J=\N^*$ sinon. 
Pour tout $j$ dans $J$ la famille admissible simple ${\mathcal S}^{(j)}$ est constitu\'ee d'une \emph{partie discr\`ete} ${\mathcal D}^{(j)}$ et d'une \emph{partie continue} ${\mathcal C}^{(j)}$,  
$${\mathcal S}^{(j)} \ = \ \bigl ( {\mathcal D}^{(j)};{\mathcal C}^{(j)} \bigr ) \ = \ \bigl ( {(\mu _l^{(j)})}_{l \in L^{(j)}} ; {(\mu _{\alpha}^{(j)})}_{\alpha \in A^{(j)}} \bigr ) \ , $$   
avec $L^{(j)}=\{ 1, \ldots , n_j\}$ ou $L^{(j)}=\N^*$ et $A^{(j)}$ ensemble totalement ordonn\'e sans \'el\'ement maximal, v\'erifiant: 

- pour $j$ appartenant \`a $J^*$, la partie discr\`ete ${\mathcal D}^{(j)} ={ \bigl (\mu _l^{(j)}\bigr )}_{l \in L^{(j)}}$ est finie et la partie continue ${\mathcal C}^{(j)} = {\bigl (\mu _{\alpha}^{(j)}\bigr )}_{\alpha \in A^{(j)}}$ est non vide, et la premi\`ere valuation $\mu _1^{(j+1)}$ de la famille simple ${\mathcal S}^{(j+1)}$ est une valuation augment\'ee limite pour la famille admissible continue ${\mathcal C}^{(j)}$; 

- la premi\`ere valuation $\mu ^{(1)}_1$ de la famille est la valuation associ\'ee \`a un polyn\^ome unitaire $\phi ^{(1)}_1$ de degr\'e un et \`a une valeur $\gamma ^{(1)}_1$, $\mu ^{(1)}_1 = [ \nu \ ; \ \mu ^{(1)}_1 ( \phi ^{(1)}_1) = \gamma ^{(1)}_1 ]$. 

Dans la suite, comme la valuation $\nu$ de $K$ est fix\'ee nous dirons simplement que ${\mathcal A}$ est une famille admissible de valuations de $K[x]$. 
\end{definition}  

La partie discr\`ete ${\mathcal D} = \bigl (\mu _i\bigr )_{i \in I}$ d'une famille admissible simple $\mathcal S$ est une famille non vide de valuations augment\'ees it\'er\'ees de $K[x]$ telle que la famille de polyn\^omes-cl\'es $\bigl (\phi _i\bigr )_{i \in I}$ v\'erifie l'in\'egalit\'e stricte deg$\phi _i>$ deg$\phi _{i-1}$.

La partie continue ${\mathcal C} = \bigl (\mu _{\alpha}\bigr )_{\alpha \in A}$ de la famille admissible simple $\mathcal S$ est une famille de valuations de $K[x]$, o\`u l'ensemble $A$ est un ensemble totalement ordonn\'e, sans \'el\'ement maximal, associ\'ee \`a une famille de polyn\^omes $\bigl (\phi _{\alpha}\bigr )_{\alpha \in A}$ de m\^eme degr\'e $d$, et \`a une famille $\bigl (\gamma _{\alpha}\bigr )_{\alpha \in A}$ d'\'el\'ements de $\tilde\Gamma$,     
pour tout $\alpha <\beta$ dans $A$, $\phi _{\beta}$ est un polyn\^ome-cl\'e pour la valuation $\mu _{\alpha}$ et la valuation $\mu _{\beta}$ est la valuation augment\'ee $\mu _{\beta}=[\mu _{\alpha} \ ; \ \mu _{\beta} (\phi _{\beta})=\gamma _{\beta}]$, avec $\gamma _{\beta}>\mu _{\alpha}(\phi _{\beta})=\gamma _{\alpha}$.  
La famille ${\mathcal C}$ est vide si la famille ${\mathcal D}$ est infinie, sinon le degr\'e $d$ des polyn\^omes-cl\'e $\phi _{\alpha}$ est \'egal au degr\'e du dernier polyn\^ome-cl\'e $\phi _n$ de la famille $\bigl (\phi _i \bigr ) _{i \in I}$ associ\'ee \`a ${\mathcal D}$, et pour tout $\alpha$ dans $A$, la valuation $\mu _{\alpha}$ est la valuation augment\'ee $\mu _{\alpha} = [\mu _n \ ; \ \mu _{\alpha} ( \phi _{\alpha} ) = \gamma _{\alpha} ]$. 

\vskip .2cm

Nous pouvons aussi \'ecrire la famille admissible ${\mathcal A}$ comme une famille index\'ee par un ensemble totalement ordonn\'e $I$,    
$${\mathcal A} \ = \ ( \mu _i )_{i \in I} \ , $$  
et l'ensemble $I$ peut \^etre d\'ecrit de la mani\`ere suivante:
pour tout $j$ dans $J$, nous munissons l'ensemble $B^{(j)} = L^{(j)} \sqcup A^{(j)}$ de l'ordre total induit par les ordres sur $L^{(j)}$ et sur $A^{(j)}$ et d\'efini par $l < \alpha$ pour tout $l \in L^{(j)}$ et tout $\alpha \in A^{(j)}$; et nous posons     
$$I \ = \ \bigl\{ (j,b) \ | \ j\in J \ \hbox{et} \ b \in B^{(j)} \bigr\}\ , $$   
muni de l'ordre lexicographique, c'est-\`a-dire $(j,b) < (j',b')$ si $j<j'$ dans $J$ et $(j,b) \leq (j,b')$ si et seulement si $b \leq b'$ dans $B^{(j)}$. 
L'ordre sur l'ensemble $I$ peut \^etre caract\'eris\'e par la relation suivante: $i<k$ dans $I$ si et seulement si pour tout polyn\^ome $f$ de $K[x]$ nous avons $\mu _i (f) \leq \mu _k(f)$ et il existe au moins un polyn\^ome $g$ avec $\mu _i(g)<\mu _k(g)$.       

\vskip .2cm

A toute famille admissible ${\mathcal A}$ nous associons la famille des polyn\^omes-cl\'es ou polyn\^omes-cl\'es limites $\bigl ( \phi _i \bigr )_{i \in I}$, que nous appelons pour simplifier la famille des polyn\^omes-cl\'es, et la famille des valeurs $\bigl ( \gamma _i \bigr )_{i \in I}$. 

\vskip .2cm

\begin{definition} 
 Une famille admissible $\mathcal{A} = ( \mu _i )_{i \in I}$ est une famille \emph{admise} si pour tout polyn\^ome $f$ dans $K[x]$ la famille $( \mu _i(f)) _{i \in I}$ admet un plus grand \'el\'ement dans le groupe $\Gamma$.
\end{definition}

Si l'ensemble $I$ poss\`ede un plus grand \'el\'ement $\bar\iota$, c'est le cas quand la famille admissible ${\mathcal A}$ est r\'eunion d'un nombre fini de familles simples et quand la derni\`ere famille simple ${\mathcal S}^{(N)}$ est discr\`ete finie,  ${\mathcal S}^{(N)} =\bigl (\mu _1^{(N)} ,\ldots ,\mu _{n_N}^{(N)} \bigr )$, alors la famille $\mathcal{A}$ est admise et est dite \emph{compl\`ete}.              
Dans ce cas la derni\`ere valuation $\mu _{\bar\iota} =  \mu _{n_N}^{(N)}$ peut \^etre une pseudo-valuation de $K[x]$. 

Si l'ensemble $I$ n'a pas de plus grand \'el\'ement, une famille $\mathcal{A} = ( \mu _i )_{i \in I}$ est admise si pour tout polyn\^ome $f$ il existe $i \in I$ tel que $\mu _i(f) = \mu _j(f)$ pour tout $j \geq i$.  
C'est le cas si la famille $\mathcal{A}$ est r\'eunion infinie de familles admissibles simples, ou si la famille $\mathcal{A}$ est r\'eunion de $t$ familles admissibles simples ${\mathcal A} = {\mathcal S}^{(1)} \cup \ldots \cup {\mathcal S}^{(N)}$,  telle que la derni\`ere famille simple ${\mathcal S}^{(N)}$ est une famille discr\`ete infinie, c'est-\`a-dire ${\mathcal S}^{(N)} = {\bigl ( \mu _l^{(N)}\bigr )}_{l \in L^{(N)}}$ avec $L^{(N)}$ infini, ou enfin si la famille simple ${\mathcal S}^{(N)}$ est une famille infinie, ${\mathcal S}^{(N)} = \bigl ( {(\mu _l^{(N)})}_{l \in L^{(N)}}; {(\mu _{\alpha}^{(N)})}_{\alpha \in A^{(N)}}\bigr )$, telle que pour tout $f$ dans $K[x]$ il existe $\alpha < \beta$ dans $A^{(N)}$ avec $\mu _{\alpha} ^{(N)}(f) = \mu  _{\beta} ^{(N)}(f)$.  

La famille admise ${\mathcal A} = ( \mu _i )_{i \in I}$ converge vers la valuation ou pseudo-valuation $\mu$ de $K[x]$ d\'efinie pour tout polyn\^ome $f$ par 
$$\mu (f) \ = \ Sup \bigl ( \mu _i (f) \ ; \ i \in I \bigr ) \ .$$ 
Si l'ensemble $I$ admet un plus grand \'el\'ement $\bar \iota$, la limite de la famille est la valuation ou pseudo-valuation $\mu _{\bar \iota}$, sinon la limie est une valuation d\'efinie par $\mu (f) = \mu _i(f)$ pour $i$ assez grand dans $I$. 

Nous avons une r\'eciproque au r\'esultat pr\'ec\'edent. 

\begin{theorem}\label{th1} (Th\'eor\`emes 2.4. et 2.5. de \cite{Va 1})
Soit $\mu$ une valuation ou pseudo-valuation de $K[x]$ prolongeant une valuation $\nu$ de $K$, alors il existe une famille admise de valuations de $K[x]$, not\'ee $\mathcal{A}(\mu )$ et appel\'ee \emph{famille admise associ\'ee \`a la valuation $\mu$} qui converge vers $\mu$. 
\end{theorem}

\begin{remark}\label{rmq:unicite}
La famille admise associ\'ee \`a une valuation $\mu$ n'est pas unique, mais est d\'etermin\'ee \`a \'equivalence pr\`es, o\`u deux familles admissibles ${\mathcal A} = \bigcup _{j \in J} {\mathcal S}^{(j)}$ et ${\mathcal A}' = \bigcup _{j \in J'} {\mathcal S}'^{(j)}$ sont dites \'equivalentes si $J=J'$, si les familles discr\`etes ${\mathcal D}^{(j)} = \bigl ( \mu ^{(j)}_i \bigr ) _{1 \leq i \leq n}$ et ${{\mathcal D}' }^{(j)} = \bigl ( {\mu '}^{(j)} _i \bigr ) _{1 \leq i \leq n'}$ co\"\i ncident jusqu'\`a l'avant-derni\`ere valuation, c'est-\`a-dire quand $n=n'$ et $\mu ^{(j)}_i = \mu '^{(j)}_i$ pour tout $i$, $1 \leq i \leq n-1$, et si les sous-familles continues ${\mathcal C}^{(j)}$ et ${\mathcal C}'^{(j)}$ co\"\i ncident asymptotiquement.  (cf. Proposition 2.9. de \cite{Va 2})
\end{remark}   

\vskip .2cm 

      \section{Ind\'ependance du g\'en\'erateur}    

Soient $L$ une extension monog\`ene de $K$ et $\tilde\mu$ une valuation de $L$ qui prolonge $\nu$. 
Si nous choisissons un g\'en\'erateur de l'extension $L$, nous pouvons \'ecrire le corps $L$ soit comme le corps des fractions de l'anneau des polyn\^omes $K[x]$, soit comme le corps quotient de l'anneau des polyn\^omes $K[x]$ par un id\'eal premier $\mathfrak{p}$. 
Dans le premier cas, $L$ est une extension transcendante pure, $L=K(x)$, et la valuation $\tilde\mu$ de $L$ induit par restriction une valuation $\mu$ de $K[x]$, dans le deuxi\`eme cas $L$ est une extension finie, $L=K(\theta )$, et la valuation $\tilde\mu$ de $L$ induit une pseudo-valuation $\mu$ de $K[x]$ dont le socle est l'id\'eal premier $\mathfrak{p}$ engendr\'e par le polyn\^ome minimal de $\theta$ sur $K$. 

La construction de la famille admise $\mathcal{A}(\mu )$ associ\'ee par le th\'eor\`eme \ref{th1} \`a la valuation ou pseudo-valuation $\mu$ de l'anneau $K[x]$ se fait par approximations successives. 
Nous supposons que nous avons construit une famille admissible $\mathcal{A}' = \bigl ( \mu _i \bigr )_{i \in I'}$ convergeant vers une valuation $\mu '$ telle que pour tout polyn\^ome $f$ de $K[x]$ et tout $i\in I'$ nous ayons l'in\'egalit\'e $\mu '(f) \leq \mu (f)$. Si la valuation $\mu '$ est diff\'erente de la valuation $\mu$ nous consid\'erons alors l'ensemble suivant 
$$\Phi _{\mu}(\mu ') \ = \ \{ \phi \in K[x] \ | \ \mu ' ( \phi ) < \mu (\phi ) \ , \ {\rm deg}\phi = d \ \hbox{et $\phi$ unitaire} \} \ ,$$ 
o\`u $d$ est le degr\'e minimal des polyn\^omes v\'erifiant l'in\'egalit\'e stricte $\mu '(f) < \mu (f)$.

Nous choisissons un un polyn\^ome $\phi _j$ de cet ensemble, nous v\'erifions que c'est un polyn\^ome-cl\'e ou un polyn\^ome-cl\'e limite qui permet de construire une valuation augment\'ee ou une valuation augment\'ee limite $\mu _j$. 
Nous obtenons ainsi une nouvelle famille admissible $\mathcal{A}''= \bigl ( \mu _i \bigr )_{i \in I''}$, o\`u l'ensemble ordonn\'e $I''$ est \'egal \`a $I '\cup \{j\}$ avec $j>i$ pour tout $i \in I'$. 
 
La nouvelle famille obtenue converge par construction vers la valuation $\mu _j$, qui est une \emph{meilleure approximation} de la vauation $\mu$ que la valuation $\mu '$ dans la mesure o\`u l'ensemble $\Phi _{\mu}(\mu _j)$ des polyn\^omes $f$ de $K[x]$ qui v\'erifient l'in\'egalit\'e stricte $\mu _j(f) < \mu (f)$ est strictement inclus dans l'ensemble $\Phi _{\mu}(\mu ')$ des polyn\^omes $f$ qui v\'erifient $\mu '(f) < \mu (f)$. 
Nous renvoyons aux articles de l'auteur \cite{Va 1} et \cite{Va 2} pour la construction pr\'ecise de la famille admise associ\'ee. 

Cette construction n'a de sens que sur l'anneau des polyn\^omes, en particulier elle fait intervenir la notion de degr\'e, la division euclidienne et une relation de comparaison de la forme $\mu _i(f) \leq \mu _j(f)$ pour tout $f$ dans $K[x]$, alors qu'il est impossible d'avoir une relation de comparaison de la forme $\mu _i(u) \leq \mu _j(u)$ pour tout \'el\'ement $u$ d'un corps avec deux valuations distinctes. 
Il semble donc que la famille admise $\mathcal{A}(\mu )$ d\'epende de mani\`ere essentielle du choix du g\'en\'erateur, et dans le d\'ebut du paragraphe 2 de \cite{Va 2} l'auteur fait la remarque suivante 
\emph{``... la famille admissible de valuations que nous allons d\'efinir \`a partir de la valuation $\mu$ de $L$ d\'epend du g\'en\'erateur $x$ ou $\theta$ de $L$ sur $K$."}  

\vskip .2cm 

En fait nous allons montrer que dans le cas o\`u $L$ est une extension transcendante pure de $K$, c'est-\`a-dire dans le cas o\`u $L$ est le corps des fractions de l'anneau des polyn\^omes $K[x]$, la famille de valuations $\mathcal{A} (\mu )= \bigl ( \mu _i \bigr )_{i \in I}$, consid\'er\'ees comme des valuations du corps $L$, est ind\'ependante du g\'en\'erateur $x$ choisi.  

Plus pr\'ecis\'ement nous allons comparer les familles admises de valuations de $L$, ${\mathcal A}_{[x]} = \bigl ( \mu _{[x],i} \bigr )_{i \in I_{[x]}}$ et ${\mathcal A}_{[y]}= \bigl ( \mu _{[y],i} \bigr )_{i \in I_{[y]}}$, associ\'ees respectivement 
aux restrictions de la valuation $\mu$ aux anneaux de polyn\^omes $K[x]$ et $K[y]$, o\`u $x$ et $y$ sont deux g\'en\'erateurs de l'extension $L$. 

Pour cela nous devons d'abord conna\^\i tre les relations qui existent entre deux g\'en\'erateurs, ce qui est donn\'e par le r\'esultat classique suivant. 

\begin{proposition}\label{prop:autL/K}    
Le groupe $Aut_K(L)$ des automorphismes de $L=K(x)$ sur $K$ est isomorphe au groupe projectif lin\'eaire $PGL_2(K)$. 

A toute matrice $M= \left(
\begin{array}{cc}
a & b \\
c & d
\end{array}
\right)$ 
est associ\'e l'automorphisme $\sigma _M$ d\'efini par 
$$\xymatrix @R=.1cm @C=1cm  {
\sigma _M : x \ar@{|->}[r] & \displaystyle{\frac{ax+b}{cx+d}}}$$

En particulier $Aut_K(L)$ est engendr\'e par les automorphismes 
$$\xymatrix @R=.1cm @C=.1cm  {
inv: x \ar@{|->}[rrrr] &&&& 1/x \\ 
t_{a,b}: x \ar@{|->}[rrrr] &&&& ax+b ,  & a \in K ^* ,\ b \in K \ .}$$
\end{proposition}   

\vskip .2cm 

\begin{theorem} 
La famille admise associ\'ee \`a la valuation $\mu$ ne d\'epend pas du g\'en\'erateur $x$ de $L$ choisi. 
\end{theorem}

\begin{preuve}  
Soient $L$ une extension transcendante pure de degr\'e un de $K$ et $\mu$ une valuation de $L$ qui prolonge la valuation $\nu$ de $K$. 
Soient $x$ et $y$ deux g\'en\'erateurs de $L$ sur $K$, nous appelons respectivement ${\mathcal A}_{[x]} = \bigl ( \mu _{[x],i} \bigr )_{i \in I_{[x]}}$ et ${\mathcal A}_{[y]}= \bigl ( \mu _{[y],i} \bigr )_{i \in I_{[y]}}$ les familles admises de valuations de $L$ associ\'ees aux restrictions de $\mu$ aux anneaux des polyn\^ome $K[x]$ et $K[y]$, et nous devons montrer que ces familles admissibles sont \'equivalentes. 

D'apr\`es la proposition \ref{prop:autL/K}, il existe une matrice $M= \left(\begin{array}{cc} a & b \\ c & d \end{array} \right)$ telle que $y = \displaystyle{\frac{ax+b}{cx+d}}$, et il suffit de consid\'erer les cas $y=ax+b$, $a \not= 0$, et $y=1/x$. 

Nous remarquons que si $y$ est \'egal \`a $ax+b$, avec $a\not= 0$, les anneaux $K[x]$ et $K[y]$ sont \'egaux et il est facile de voir que par construction nous pouvons trouver des familles ${\mathcal A}_{[x]} = \bigl ( \mu _{[x],i} \bigr )_{i \in I_{[x]}}$ et ${\mathcal A}_{[y]}= \bigl ( \mu _{[y],i} \bigr )_{i \in I_{[y]}}$ identiques. 

\vskip .2cm 

Nous supposons maintenant que $y$ est \'egal \`a $1/x$. 
Nous appelons $IU(K[x])$, respectivement $IU(K[y])$, l'ensemble des polyn\^omes irr\'eductibles unitaires de $K[x]$, respectivement de $K[y]$, et nous d\'efinissons l'application $Inv_{x/y}$ de $IU(K[x])$ dans $K[y]$ de la mani\`ere suivante: 
$$\begin{array}{rcl} 
Inv_{x/y}(x) & = & y \\ 
Inv_{x/y}(\phi (x)) & = & \psi (y) \ , 
\end{array}$$ 
avec 
$$\psi (y)  \ = (a_0 ^{-1} y^d) \phi (x) \ = \ y^d + (a_0^{-1}a_1) y^{d -1} + \ldots + (a_0^{-1}a_{d-1}) y+ a_0^{-1}$$  
pour 
$$\phi (x) \ = \ x^d+ a_{d -1} x^{d -1} + \ldots + a_1 x + a_0 \ , \hbox{avec}\ d\geq 1 \ \hbox{et}\ a_0 \not= 0 \ .$$ 

L'application $Inv_{x/y}$ induit une bijection de $IU(K[x])$ dans $IU(K[y])$, et si nous d\'efinissons de mani\`ere analogue l'application $Inv_{y/x}$ de $IU(K[y])$ dans $IU(K[x])$, nous avons $(Inv_{x/y} )^{-1} = Inv_{y/x}$. 

\vskip .2cm 

Nous notons $\Phi _1 (K[x])$ et $\Phi _1 (K[y])$ les ensembles des polyn\^omes unitaires de degr\'e un respectivement dans $K[x]$ et dans $K[y]$, c'est-\`a-dire les ensembles des polyn\^omes de degr\'e un appartenant respectivement \`a $IU(K[x])$ et \`a $IU(K[y])$. 
Pour toute valuation $\mu$ de $K(x)$ nous d\'efinissons les deux sous-ensembles du groupe des valeurs de $\mu$, $\Lambda _{\mu ,1}(K[x]) = \{ \mu (\phi ) \ | \ \phi \in \Phi _1(K[x]) \}$ et $\Lambda _{\mu ,1}(K[y]) = \{ \mu (\psi ) \ | \ \psi \in \Phi _1(K[y]) \}$. 

\vskip .2cm 

\begin{lemma} \label{le:degre1}
Si l'ensemble $\Lambda _{\mu ,1}(K[x])$ a un plus grand \'el\'ement $\gamma _1$ et si $\phi _1$ est un polyn\^ome unitaire de degr\'e un dans $K[x]$ tel que $\mu (\phi _1) = \gamma _1$, alors l'ensemble $\Lambda _{\mu ,1}(K[y])$ a un plus grand \'el\'ement $\delta _1$ et le polyn\^ome unitaire de degr\'e un $\psi _1$ d\'efini par $\psi _1 = Inv_{x/y}( \phi _1)$ v\'erifie $\mu (\psi _1) = \delta _1$. 

Si l'ensemble $\Lambda _{\mu ,1}(K[x])$ n'a pas de plus grand \'el\'ement et si $(\phi _{\alpha})_{\alpha \in A}$ est une famille de polyn\^omes unitaires de degr\'e un dans $K[x]$ index\'ee par un ensemble totalement ordonn\'e $A$ sans plus grand \'el\'ement, telle que la famille $(\gamma _{\alpha})_{\alpha \in A}$, avec $\gamma _{\alpha} = \mu ( \phi _{\alpha})$, v\'erifie $\gamma _{\beta} > \gamma _{\alpha}$ pour $\beta > \alpha$ et est cofinale dans $\Lambda _{\mu ,1}(K[x])$, alors l'ensemble $\Lambda _{\mu ,1}(K[y])$ n'a pas de plus grand \'el\'ement et la famille $(\delta_{\alpha})_{\alpha \in A}$ d\'efinie par $\delta _{\alpha} = \mu ( \psi _{\alpha})$ o\`u $ \psi _{\alpha} = Inv_{x/y}(\phi _{\alpha})$ est cofinale dans $\Lambda _{\mu ,1}(K[y])$.   
\end{lemma} 

\begin{preuvele}
Supposons d'abord que l'ensemble $\Lambda _{\mu ,1}(K[x])$ a un plus grand \'el\'ement, c'est-\`a-dire qu'il existe un polyn\^ome unitaire de degr\'e un $\phi _1$ dans $K[x]$ tel que $\mu (\phi _1) \geq \mu (\phi )$ pour tout $\phi \in \Phi _1(K[x])$. 

Consid\'erons le cas $\phi _1 (x) =x$, c'est-\`a-dire le cas $\mu (x) \geq \mu (x +c)$ pour tout $c \in K$, et par d\'efinition le polyn\^ome $\psi _1 = Inv_{x/y}(\phi _1)$ est \'egal \`a $\psi _1 (y) =y$. 
Nous avons alors pour tout $d \in K^*$ l'in\'egalit\'e 
$$\mu ( y+d) = \mu (dy) + \mu (x+d^{-1}) \leq \mu (d) + \mu (y) + \mu (x) = \mu (d) \ , $$ 
dont nous d\'eduisons l'in\'egalit\'e 
$$\mu ( y+d) \leq \mu (y) \ ,$$ 
c'est-\`a-dire $\mu (\psi _1) \geq \mu (\psi )$ pour tout $\psi \in \Phi _1(K[y])$. 

Consid\'erons maintenant le cas $\phi _1 (x) =x + c_0$ avec $c_0 \not= 0$, nous avons $\mu (x+c_0) \geq \mu (x +c)$ pour tout $c \in K$, et nous pouvons supposer que nous avons $\mu (x+c_0) > \mu (x)$, d'o\`u $\mu (x) = \mu (c_0)$. 
Le polyn\^ome $\psi _1  = Inv_{x/y}(\phi _1)$ est alors \'egal \`a $\psi _1(y) = y + d_0$ avec $d_0 = {c_0}^{-1}$, et nous avons $\mu (y) = \mu (d_0) = - \mu (x)$. 
 
Supposons qu'il existe $d \in K$ tel que $\mu (y+d) > \mu (y+ d_0)$, d'o\`u $\mu (y+d_0) = \mu (d-d_0)$. 
Nous en d\'eduisons 
$$\begin{array}{rcl}
\mu (x+c_0) = \mu (xc_0) + \mu (y+ d_0) & < & \mu (xc_0) + \mu (y+ d) \\ 
& = & \mu (x c_0) + \mu (yd) + \mu (x + d^{-1}) \\
& \leq & \mu (c_0d) + \mu (x + c_0) \ , 
\end{array}$$ 
d'o\`u $\mu (c_0d) = \mu (d/d_0) >0$ et $\mu (d) > \mu (d_0) = \mu ( d-d_0)$.   
Nous avons alors les deux in\'egalit\'es $\mu (d) > \mu (y)$ et $\mu (y+d) > \mu (y)$, ce qui est impossible. 
Par cons\'equent le polyn\^ome $\psi _1$ v\'erifie bien $\mu (\psi _1) \geq \mu (\psi )$ pour tout $\psi \in \Phi _1(K[y])$. 

\vskip .2cm 

Supposons maintenant que l'ensemble $\Lambda _{\mu ,1}(K[x])$ n'a pas de plus grand \'el\'ement et soit $(\phi _{\alpha})_{\alpha \in A}$ une famille de polyn\^omes de $K[x]$ v\'erifiant les hypoth\`eses du lemme. 
Nous pouvons \'ecrire $\phi _{\alpha} (x)= x +c _{\alpha}$, et quitte \`a se restreindre \`a un sous ensemble cofinal dans $A$ nous pouvons supposer que pour tout $\alpha \in A$ nous avons l'in\'egalit\'e $\gamma _{\alpha} = \mu ( \phi _{\alpha} ) > \mu (x)$, d'o\`u $\mu (c_{\alpha}) = \mu (x)$.  

Pour tout $\alpha$ le polyn\^ome unitaire $\psi _{\alpha}$ de $\Phi _1(K[y])$ d\'efini par $\psi _{\alpha}= Inv_{x/y}(\phi _{\alpha})$ est \'egal \`a $\psi _{\alpha} (y) = (y+ d_{\alpha})$, avec $d_{\alpha} = c _{\alpha}^{-1}$. 
En particulier nous avons $\psi _{\alpha} = \phi _{\alpha} (xc_{\alpha})^{-1}$, $\delta _{\alpha} = \mu ( \psi _{\alpha})$ est \'egal \`a $\gamma _{\alpha} - 2 \mu (x)$, et nous d\'eduisons de l'in\'egalit\'e $\gamma _{\alpha} = \mu ( \phi _{\alpha} ) > \mu (x)$ que nous avons $\delta _{\alpha} > \mu (y) = \mu ( d _{\alpha})$, par cons\'equent pour tout $\beta > \alpha$ nous avons l'in\'egalit\'e $\delta _{\beta} > \delta _{\alpha}$. 

Si $\psi (y) = y+d$ v\'erifie $\mu (\psi ) > \mu (\psi _{\alpha})$ pour $\alpha \in A$, le polyn\^ome $\psi$ est \'egal \`a $Inv_{x/y}(\phi )$ pour un polyn\^ome $\phi (x)=x+c$, avec l'\'egalit\'e $\psi = (yd) \phi$ et $c = d ^{-1}$. 
Il existe alors $\beta$ dans $A$ tel que $\gamma _{\beta} > \mu (\phi )$ d'o\`u $\delta _{\beta} > \mu (\psi )$, nous en d\'eduisons que la famille $(\delta _{\alpha})_{\alpha \in A}$ est cofinale dans $\Lambda _{\mu ,1}(K[y])$.  

\hfill \tf    
\end{preuvele}

Montrons que les premi\`eres valuations des familles ${\mathcal A}_{[x]}$ et ${\mathcal A}_{[y]}$ co\"\i ncident. 
Nous notons respectivement ${\mathcal S}^{(j)}_{[x]} = \bigl ( {\mathcal D}^{(j)}_{[x]};{\mathcal C}^{(j)}_{[x]} \bigr ) $ et ${\mathcal S}^{(j)}_{[y]} = \bigl ( {\mathcal D}^{(j)}_{[y]};{\mathcal C}^{(j)}_{[y]} \bigr ) $ les sous-familles admissibles simples de ${\mathcal A}_{[x]}$ et ${\mathcal A}_{[y]}$, constitu\'ees de leur parties discr\`etes et continues. 

Nous rappelons que pour tout polyn\^ome $\phi$ dans $\Phi _1(K[x])$ et toute valeur $\gamma$ nous pouvons d\'efinir une valuation $\mu _{[x]}^{[\phi ,\gamma ]}$ de $L$ comme la valuation augment\'ee de $K[x]$ d\'efinie par 
$$\mu _{[x]}^{[\phi ,\gamma ]} \ = \ [ \nu \ ; \ \mu _{[x]}^{[\phi ,\gamma ]}(\phi ) = \gamma ] \ .$$ 
Plus pr\'ecis\'ement pour tout polyn\^ome $f$ dans $K[x]$ nous avons  
$$\mu _{[x]}^{[\phi ,\gamma ]}(f) = Inf \left ( \nu (a_j) + j \gamma  \ ; \  0 \leq j \leq n \right ) \ , $$   
o\`u $f = a_n \phi ^n + \ldots + a_1 \phi  + a_0$, et $a_j \in K$. 

La partie discr\`ete $\mathcal{D}^{(1)}_{[x]}$ de la premi\`ere sous-famille admissible simple de ${\mathcal A}_{[x]}$ n'est pas r\'eduite \`a une seule valuation si l'ensemble $\Lambda _{\mu ,1}(K[x]) = \{ \mu (\phi ) \ | \ \phi \in \Phi _1(K[x]) \} $ a un plus grand \'el\'ement $\gamma _1$. Dans ce cas la premi\`ere valuation $\mu _{[x],1}^{(1)}$ de la famille ${\mathcal A}_{[x]}$ est la premi\`ere valuation de $\mathcal{D}^{(1)}_{[x]}$ et est \'egale \`a 
$$ \mu _{[x],1}^{(1)} \ = \ \mu _{[x]}^{[\phi _1,\gamma _1]} \ , $$ 
o\`u $\phi _1$ est un polyn\^ome de $\Phi _1(K[x])$ tel que $\mu (\phi _1) = \gamma _1$. 

De plus si la valuation $\mu$ n'est pas \'egale \`a $\mu _{[x],1}^{(1)}$, la deuxi\`eme valuation $\mu _{[x],2}^{(1)}$ de la famille ${\mathcal A}_{[x]}$ est obtenue comme valuation augment\'ee de la valuation $\mu _{[x],1}^{(1)}$ associ\'ee \`a un polyn\^ome $\phi _2$ de degr\'e $d \geq 2$: 
$$\mu _{[x],2}^{(1)} \ = \ [ \mu _{[x],1}^{(1)}  \ ; \ \mu _{[x],2}^{(1)}  (\phi _2) = \gamma _2 ] \ .$$ 

La partie discr\`ete $\mathcal{D}^{(1)}_{[x]}$ est r\'eduite \`a une seule valuation si l'ensemble $\Lambda _{\mu ,1}(K[x])$ n'a pas de plus grand \'el\'ement. Dans ce cas nous pouvons choisir comme premi\`ere valuation $\mu _{[x],1}^{(1)}$ de la famille $\mathcal{A}_{[x]}$ la valuation augment\'ee associ\'ee \`a n'importe quel polyn\^ome $\phi _1$ de $\Phi _1(K[x])$ et \`a la valeur $\gamma _1 = \mu (\phi _1)$:  
$$\mu _{[x],1}^{(1)} \ = \ [ \nu \ ; \ \mu _{[x],1}^{(1)}(\phi _1) = \gamma _1] \ ,$$ 
ainsi qu'une famille $(\phi _{\alpha})_{\alpha \in A}$ de polyn\^omes unitaires de degr\'e un dans $K[x]$ index\'ee par un ensemble totalement ordonn\'e $A$ sans plus grand \'el\'ement, telle que la famille $(\gamma _{\alpha})_{\alpha \in A}$, avec $\gamma _{\alpha} = \mu ( \phi _{\alpha})$, v\'erifie $\gamma _{\beta} > \gamma _{\alpha} > \gamma _1$ pour tout $\beta > \alpha$ et soit cofinale dans $\Lambda _{\mu ,1}(K[x])$. 
La partie continue $\mathcal{C}^{(1)}_{[x]}$ est alors constitu\'ee des valuations augment\'ee $\mu _{[x],\alpha}^{(1)}$ associ\'ees aux polyn\^omes  $\phi _{\alpha}$  et aux valeurs $\gamma _{\alpha}$ :
$$ \mu _{[x],\alpha}^{(1)} \ = \ [ \nu  \ ; \ \mu _{[x],\alpha} ^{(1)} (\phi _{\alpha}) = \gamma _{\alpha} ] \ = \ [ \mu _{[x],1}^{(1)}  \ ; \ \mu _{[x],\alpha} ^{(1)} (\phi _{\alpha}) = \gamma _{\alpha} ]\ .$$ 

De plus si la valuation $\mu$ n'est pas \'egale \`a la valuation limite, c'est-\`a-dire s'il existe $f \in K[x]$ tel que $\mu _{[x],\alpha} ^{(1)}(f) < \mu (f)$ pour tout $\alpha$ dans $A$, la famille $\mathcal{A}_{[x]}$ contient une deuxi\`eme sous-famille admissible simple et la premi\`eme valuation $\mu _{[x],1}^{(2)}$ de la partie discr\`ete $\mathcal{D}^{(2)}_{[x]}$ est obtenue comme valuation augment\'ee limite de la famille $\bigl (\mu _{[x],\alpha}^{(1)} \bigr )_{\alpha \in A}$ associ\'ee \`a un polyn\^ome $\phi _1^{(2)}$ de degr\'e $d \geq 2$: 
$$\mu _{[x],1}^{(2)}\ = \ \Bigl [ \bigl (\mu _{[x],\alpha}^{(1)} \bigr )_{\alpha \in A}  \ ; \ \mu _{[x],1}^{(2)}  (\phi _1^{(2)}) = \gamma _1^{(2)} \Bigr ] \ .$$ 

\vskip .2cm 

Supposons d'abord que l'ensemble $\Lambda _{\mu ,1}(K[x]) = \{ \mu (\phi ) \ | \ \phi \in \Phi _1(K[x]) \} $ a un plus grand \'el\'ement $\gamma _1$, et soit $\phi _1$ un polyn\^ome de $\Phi _1(K[x])$ tel que $\mu (\phi _1) = \gamma _1$, alors la valuation $\mu _{[x],1}^{(1)}$ est la valuation augment\'ee $ \mu _{[x],1}^{(1)} \ = \ \mu _{[x]}^{[\phi _1,\gamma _1]}$ de $K[x]$.   

D'apr\`es le lemme \ref{le:degre1} l'ensemble $\Lambda _{\mu ,1}(K[y]) = \{ \mu (\psi ) \ | \ \psi \in \Phi _1(K[y]) \} $ a un plus grand \'el\'ement $\delta _1$, et le polyn\^ome unitaire de degr\'e un $\psi _1$ dans $K[y]$ d\'efini par $\psi _1 = Inv_{x/y}( \phi _1)$ v\'erifie $\mu (\psi _1) = \delta _1$. 
Alors de la m\^eme mani\`ere la valuation $\mu _{[y],1}^{(1)}$ est la valuation augment\'ee $ \mu _{[y],1}^{(1)} \ = \ \mu _{[y]}^{[\psi _1,\delta _1]}$ de $K[y]$. 

\vskip .2cm

Si le polyn\^ome $\phi _1(x)$ est \'egal \`a $x$ alors $\psi _1$ est \'egal \`a y et $\mu (y) =\delta _1 = - \gamma _1 = - \mu (x)$. 
Pour tout polyn\^ome $f = a_n x ^n + \ldots + a_1 x + a_0$ dans $K[x]$, nous posons $g = y^n f$, d'o\`u $g = a_0 y ^n + \ldots + a_{n-1} y + a_n$, et nous trouvons d'apr\`es ce qui pr\'ec\`ede 
$$\begin{array}{rcl}
\mu _{[y],1}^{(1)}(f) & = & \mu  _{[y],1}^{(1)}(y^{-n} g) \\ 
& = & - n \delta _1 + Inf \left ( \nu (a_j) + (n-j) \delta _1 \ ; \  0 \leq j \leq n \right ) \\
 & = & Inf \left ( \nu (a_j) + j \gamma _1 \ ; \  0 \leq j \leq n \right )   \\   
 & = & \mu _{[x],1}^{(1)}(f) \ . 
\end{array}$$

\vskip .2cm

Si le polyn\^ome $\phi _1(x)$ est \'egal \`a $x+c_0$ avec $c_0 \not=0$, alors $\psi _1(y)$ est \'egal \`a $y+d_0$ avec $d_0 = c_0^{-1}$ si et $\psi _1 = (yd_0) \phi _1$, et nous pouvons supposer $\mu (\phi _1) > \mu (x)= \mu (c_0)$, d'o\`u $\mu (\psi _1) > \mu (y) = \mu (d_0)$. 
Si nous posons $\gamma _0 = \mu (x) = - \mu (y) = - \delta _0$ nous avons l'\'egalit\'e 
$$\mu ( \psi _1)  = \delta _1 = \gamma _1 + 2 \delta _0 = \mu (\phi _1) +2 \delta _0 \ . $$

Soit $\mu _{[x],0}^{(1)}$ la valuation de $K[x]$ d\'efinie par   
$$\mu _{[x],0}^{(1)} \ = \ [ \nu \ ; \ \mu _{[x],0}^{(1)}(\phi _0) = \gamma _0] \ ,$$  
avec $\phi _0 =x$ et $\gamma _0 = \mu (x)$, alors le polyn\^ome $\phi _1$ est un polyn\^ome-cl\'e pour $\mu _{[x],0}^{(1)}$ et $\mu _{[x],0}^{(1)} (\phi _1) = \gamma _0$. 
D'apr\`es le lemme 15.1. de \cite {McL 1} la valuation $\mu _{[x],1}^{(1)}$ est \'egale \`a la valuation augment\'ee $[ \mu _{[x],0}^{(1)} \ ; \ \mu _{[x],1}^{(1)}(\phi _1) = \gamma _1] $ et comme $x$ n'est pas $\mu _{[x],0}^{(1)}$-divisible par $\phi _1$, nous d\'eduisons du th\'eor\`eme 5.2. de \cite {McL 1} que pour tout polyn\^ome $f$ de $K[x]$ qui s'\'ecrit sous la forme $f = a_n x ^{e_n} \phi _1 ^n + \ldots + a_1 x^{e_1} \phi _1 + a_0 x^{e_0}$, avec $e_j \in \N$ et $a_j \in K$ nous avons l'\'egalit\'e 
$$\mu _{[x],1}^{(1)}(f) = Inf \left ( \nu (a_j) + e_j \gamma _0 + j \gamma  _1\ ; \  0 \leq j \leq n \right ) \ . $$   
De m\^eme, si nous posons $\mu (y) = \delta _0$,  pour tout polyn\^ome $g$ de $K[y]$ qui s'\'ecrit sous la forme $g = b_n y ^{e_n} \psi _1 ^n + \ldots + b_1 y^{e_1} \psi _1 + b_0 y^{e_0}$, avec $e_j \in \N$ et $b_j \in K$ nous avons l'\'egalit\'e 
$$\mu _{[y],1}^{(1)}(g) = Inf \left ( \nu (b_j) + e_j \delta _0 + j \delta  _1\ ; \  0 \leq j \leq n \right ) \ . $$   

Pour tout polyn\^ome $f = a_n \phi _1 ^n + \ldots + a_1 \phi _1 + a_0$ dans $K[x]$, avec $a_j \in K$, nous posons $g = (yd_0)^n f$, d'o\`u 
$$\begin{array}{rcl} 
g & = & a_0 (yd_0 \phi _1) ^n + \ldots + a_j (yd_0)^{n-j} ( yd_0 \phi _1) ^j + \dots + a_{n-1} (yd_0)^{n-1} ( yd_0 \phi _1) + a_n (yd_0)^n  \\ 
& = & a_0 \psi _1 ^n + \ldots + a_j (yd_0)^{n-j} \psi _1 ^j + \dots + a_{n-1} (yd_0)^{n-1} \psi _1 + a_n (yd_0)^n \ , 
\end{array}$$ 
et nous trouvons d'apr\`es ce qui pr\'ec\`ede 
$$\begin{array}{rcl}
\mu _{[y],1}^{(1)}(f) & = & \mu  _{[y],1}((yd_0)^{-n} g) \\ 
& = & - 2n \delta _ 0 + Inf \left ( \nu (a_j) + 2(n-j) \delta _0 + j \delta _1 \ ; \  0 \leq j \leq n \right ) \\
 & = & Inf \left ( \nu (a_j) + j (\delta _1 -2 \delta _0) \ ; \  0 \leq j \leq n \right )   \\   
  & = & Inf \left ( \nu (a_j) + j \gamma _1 \ ; \  0 \leq j \leq n \right )   \\  
 & = & \mu _{[x],1}^{(1)}(f) \ . 
\end{array}$$

\vskip .2cm

Supposons qu'il n'existe pas de polyn\^ome unitaire de degr\'e un $\phi _1$ dans $K[x]$ tel que $\mu (\phi _1) \geq \mu (x+c)$ pour tout $c \in K$ et soit comme pr\'ec\'edemment $(\phi _{\alpha})_{\alpha \in A}$ une famille de polyn\^omes unitaires de degr\'e un dans $K[x]$ index\'ee par un ensemble totalement ordonn\'e $A$ sans plus grand \'el\'ement, telle que la famille $(\gamma _{\alpha})_{\alpha \in A}$, avec $\gamma _{\alpha} = \mu ( \phi _{\alpha})$, v\'erifie $\gamma _{\beta} > \gamma _{\alpha}$ pour $\beta > \alpha$ et soit cofinale dans $\Lambda _{\mu ,1}(K[x]) = \{ \mu (\phi ) \ | \ \phi \in \Phi _1(K[x]) \} $, et v\'erifie $\mu (\phi _{\alpha}) > \mu (x)$ pour tout $\alpha$ dans $A$.  

Posons $A_{[x]}^{(1)} = A$, alors pour tout $\alpha$ dans $A_{[x]}^{(1)}$ nous d\'efinissons la valuation augment\'ee $\mu _{[x],\alpha}^{(1)} = [ \nu \ ; \ \mu _{[x],\alpha}^{(1)}(\phi _{\alpha}) = \gamma _{\alpha}]$ sur $K[x]$, et la famille $\mathcal{C}^{(1)}_{[x]} = \bigl ( \mu _{[x],\alpha}^{(1)} \bigr ) _{\alpha \in A_{[x]}^{(1)}}$ est une famille admissible continue. 
Nous pouvons choisir comme premi\`ere famille admissible simple $\mathcal{S}^{(1)}_{[x]}$ de la famille admise $\mathcal{A}_{[x]}$ la r\'eunion 
$${\mathcal S}^{(1)}_{[x]} \ = \ \left ( {\mathcal D}^{(1)}_{[x]} ;{\mathcal C}^{(1)}_{[x]} \right ) \ = \left (  \bigl ( {\mu _{[x],1}^{(1)} \bigr ) ;  \bigl ( \mu _{[x],\alpha}^{(1)}} \bigr )_{\alpha \in A^{(1)}_{[x]}} \right ) \ , $$   
o\`u la valuation augment\'ee $\mu _{[x],1}^{(1)}$ est la valuation associ\'ee au polyn\^ome $\phi _{[x],1}^{(1)} = x$. 

D'apr\`es le lemme pr\'ec\'edent l'ensemble $\Lambda _{\mu ,1}(K[y])$ n'a pas de plus grand \'el\'ement et la famille $(\delta_{\alpha})_{\alpha \in A}$ d\'efinie par $\psi _{\alpha} = \mu \bigl (Inv_{x/y}(\phi _{\alpha}) \bigr )$ est cofinale dans $\Lambda _{\mu ,1}(K[y])$, et v\'erifie $\mu (\psi _{\alpha}) > \mu (y)$ pour tout $\alpha$ dans $A$.  

De mani\`ere analogue nous posons $A_{[y]}^{(1)} = A$, pour tout $\alpha$ dans $A_{[y]}^{(1)}$ nous d\'efinissons la valuation augment\'ee $\mu _{[y],\alpha}^{(1)} = [ \nu \ ; \ \mu _{[y],\alpha}^{(1)}(\psi _{\alpha}) = \delta _{\alpha}]$ sur $K[y]$, et nous pouvons choisir comme premi\`ere famille admissible simple $\mathcal{S}^{(1)}_{[y]}$ de la famille admise $\mathcal{A}_{[y]}$ la r\'eunion 
$${\mathcal S}^{(1)}_{[y]} \ = \ \left ( {\mathcal D}^{(1)}_{[y]} ;{\mathcal C}^{(1)}_{[y]} \right ) \ = \left (  \bigl ( {\mu _{[y],1}^{(1)} \bigr ) ;  \bigl ( \mu _{[y],\alpha}^{(1)}} \bigr )_{\alpha \in A^{(1)}_{[y]}} \right ) \ , $$   
o\`u la valuation augment\'ee $\mu _{[y],1}^{(1)}$ est la valuation associ\'ee au polyn\^ome $\psi _{[y],1}^{(1)} = y$. 

Comme nous avons suppos\'e $\mu (\phi _{\alpha}) > \mu (x)$ pour tout $\alpha$ dans $A$, nous montrons comme pr\'ec\'edemment que les valuations $\mu _{[x], \alpha}^{(1)}$ et $\mu _{[y], \alpha}^{(1)}$ sont \'egales comme valuations sur $L$. 
Par cons\'equent les deux premi\`eres familles admissibles simples ${\mathcal S}^{(1)}_{[x]}$ et ${\mathcal S}^{(1)}_{[y]}$ sont \'egales. 

\vskip .2cm 

Supposons maintenant que les familles ${\mathcal A}_{[x]} = \bigl ( \mu _{[x],i} \bigr )_{i \in I_{[x]}}$ et ${\mathcal A}_{[y]}= \bigl ( \mu _{[y],i} \bigr )_{i \in I_{[y]}}$ sont les m\^emes jusqu'\`a l'indice $k$, plus pr\'ecis\'ement les ensembles $\{ j \in I_{[x]}\ | \ j \leq k\}$ et $\{ j \in I_{[y]}\ | \ j \leq k\}$ sont identiques,  pour tout $j$ dans ces ensembles les valuations $\mu _{[x],j}$ et $\mu _{[y],j}$ sont \'egales comme valuations de $L$, et de plus si la valuation $\mu _{[x],j}$ est d\'efinie par le polyn\^ome-cl\'e, respectivement le polyn\^ome-cl\'e limite, $\phi _j(x)$ de degr\'e $d_j$, la valuation $\mu _{[y],j}$ est d\'efinie par le polyn\^ome-cl\'e, respectivement le polyn\^ome-cl\'e limite, $\psi _j(y)$ de degr\'e $d_j$ avec $\psi _j(y) = Inv_{x/y} (\phi _j(x))$. 
Nous notons $\mu _j$ la valuation de $L$ \'egale \`a $\mu _{[x],j}=\mu _{[y],j}$. 

De plus les valeurs associ\'ees aux valuations augment\'ees, $\gamma _j = \mu (\phi _j) = \mu _k(\phi _j)$ et $\delta _j = \mu (\psi _j) = \mu _k(\psi _j)$, sont reli\'ees par les relations suivantes:  
$$\mu (\phi _k) - \mu _j(\phi _k)  = \mu (\psi _k) - \mu _j(\psi _k) \quad \hbox{pour tout}\ 1\leq j \leq k  \ .$$

Nous rappelons la construction de la famille admise ${\mathcal A}_{[x]} = \bigl ( \mu _{[x],i} \bigr )_{i \in I_{[x]}}$ associ\'ee \`a la valuation $\mu$, si la valuation $\mu _{[x],k}$ n'est pas \'egale \`a la valuation $\mu$, l'ensemble 
$$\tilde\Phi _{\mu}(\mu _{[x],k}) \ = \ \{ f \in K[x] \ | \ \mu _{[x],k} (f) < \mu (f) \} \ ,$$ 
est non vide et pour d\'eterminer les valuations $\mu _{[x],l}$ pour $l>k$ nous devons consid\'erer l'ensemble suivant
$$\Phi _{\mu}(\mu _{[x],k}) \ = \ \{ \phi \in K[x] \ | \ \mu _{[x],k} ( \phi ) < \mu (\phi ) \ , \ {\rm deg}\phi = d _{[x],k} \ \hbox{et $\phi$ unitaire} \} \ ,$$ 
o\`u $d _{[x],k}$ est le degr\'e minimal des polyn\^omes appartenant \`a $\tilde\Phi _{\mu}(\mu _{[x],k})$, et $d _{[x],k} \geq d_k$. 

Tout polyn\^ome $\phi$ appartenant \`a  $\Phi _{\mu}(\mu _{[x],k})$ est irr\'eductible, de degr\'e $d _{[x],k}\geq 1$, avec $\phi \not= x$, en particulier il est de la forme $\phi (x) \ = \ x^{d_{[x],k}}+ \ldots + a_1 x + a_0 $ avec $a_0 \not= 0$. 

\begin{lemma} \label{le-bij}
L'application $Inv_{x/y}$ induit une bijection de $\Phi _{\mu}(\mu _{[x],k})$ dans $\Phi _{\mu}(\mu _{[y],k})$. 
\end{lemma} 

\begin{preuvele} 
Pour tout polyn\^ome $f(x)$ dans $IU(K[x])$ de degr\'e $d$, le polyn\^ome $g(y)=Inv_{x/y}(f(x)) = (a_0^{-1} y^d) f(x)$ est de m\^eme degr\'e $d$ et nous d\'eduisons de l'\'egalit\'e 
$\mu _i (a_0^{-1} y^d) = \mu (a_0^{-1} y^d)$ 
l'\'equivalence suivante:  
$$\mu _{[x],k} (f(x)) < \mu (f(x)) \ \Longleftrightarrow \ \mu _{[y],k} (g(y)) < \mu (g(y)) \ .$$ 

\hfill \tf    
\end{preuvele} 

Soient $\phi$ un polyn\^ome appartenant \`a $\Phi _{\mu}(\mu _{[x],k})$, avec $\gamma = \mu (\phi)$, et $\psi$ son image $\psi = Inv_{x/y}(\phi)$ qui appartient \`a $\Phi _{\mu}(\mu _{[y],k})$, avec $\delta = \mu (\psi)$. 
Nous consid\'erons alors les valuations $\mu_{\phi ,\gamma}$ et $\mu_{\psi ,\delta}$ de $L$, d\'efinies respectivement comme valuation augment\'ee pour la valuation $\mu _i$ sur $K[x]$ associ\'ee au polyn\^ome-cl\'e $\phi$ et \`a la valeur $\gamma$, et sur $K[y]$ associ\'ee au polyn\^ome-cl\'e $\psi$ et \`a la valeur $\delta$: 
$$\mu_{\phi ,\gamma} = [ \mu _k\ ;\ \mu_{\phi ,\gamma} (\phi ) = \gamma ] 
\quad \hbox{et} \quad 
\mu_{\psi ,\delta} = [ \mu _k\ ;\ \mu_{\psi ,\delta} (\psi ) = \delta ] \ . $$

\begin{lemma}\label{le-valuations-egales} 
Les valuations $\mu_{\phi ,\gamma}$ et $\mu_{\psi ,\delta}$ sont \'egales. 
\end{lemma} 

\begin{preuvele} 
Nous appelons $d_{>k}$ le degr\'e des polyn\^omes $\phi$ de $K[x]$ et $\psi$ de $K[y]$, $d_{>k} = d_{[x],k} = d_{[y],k}$ , nous notons $a$ le coefficient de degr\'e $0$ de $\phi$, $a \in K^*$, de telle fa\c con que nous ayons l'\'egalit\'e 
$$\psi = Inv_{x/y}(\phi ) = (a^{-1} y^{d_{>k}}) \phi \ .$$

Soit $f$ appartenant \`a $K[x]$, et soit $f = f_m \phi ^m + \ldots + f_1 \phi + f_0$ le d\'eveloppement de $f$ selon les puissances de $\phi$. 
Alors par d\'efinition de la valuation augment\'ee $\mu _{\phi , \gamma }$ nous avons 
$$\mu_{\phi , \gamma} (f) = Inf \bigl ( \mu _i(f_j) + j \gamma \ ; \ 0 \leq j \leq m \bigr ) \ . $$

Soit $s$ le degr\'e du polyn\^ome $f_m$, par cons\'equent le degr\'e de $f$ est \'egal \`a $n=s+md_{>k}$, et l'\'el\'ement $g=y^nf$ de $L$ appartient \`a l'anneau de polyn\^omes $K[y]$. 
Et dans $K[y]$ nous avons l'\'egalit\'e suivante: 
$$g \ = \ (a^m y^s f_m) (a^{-1} y^{d_{>k}} \phi ) ^m + \ldots + (a^j y^{n-jd_{>k}} f_j) (a^{-1} y^{d_{>k}} \phi ) ^j + \ldots + (a^m y^n f_0) \ ,$$
c'est-\`a-dire que nous pouvons \'ecrire $g = g_m \psi ^m + \ldots + g_1 \phi + g_0$ o\`u chaque polyn\^ome $g_j$ de $K[y]$ v\'erifie $\mu _k(g_j) = \mu (g_j)$. 
Nous d\'eduisons du th\'eor\`eme 5.2. \cite{McL 1} que nous avons l'\'egalit\'e 
$$\mu_{\psi , \delta} (g) = Inf \bigl ( \mu _k(g_j) + j \delta \ ; \ 0 \leq j \leq m \bigr ) \ . $$
Par construction pour tout $j$ nous avons 
$$\mu _k(g_j) + j \delta \ = \ \mu _k(f_j) + j \gamma + \mu _k(y^n) \ , $$
par cons\'equent nous en d\'eduisons bien $\mu _{\phi ,\gamma}(f) = \mu_{\psi , \delta} (f)$ pour tout $f$ dans $K[x]$, donc pour tout $f$ dans $L$. 

\hfill \tf    
\end{preuvele} 

Nous revenons aux deux sous-ensembles $\Phi _{\mu}(\mu _{[x],k})$ de $K[x]$ et $\Phi _{\mu}(\mu _{[y],k})$ de $K[y]$ form\'es de polyn\^omes unitaires de m\^eme degr\'e $d_{>k}$, et aux deux sous-ensembles du groupe des valeurs $\Lambda _{\mu}(\mu _{[x],k}) = \{ \mu (\phi ) \ | \ \phi \in \Phi _{\mu}(\mu _{[x],k}) \}$ et $\Lambda _{\mu}(\mu _{[y],k}) = \{ \mu (\psi ) \ | \ \psi \in \Phi _{\mu}(\mu _{[y],k}) \}$. 

\vskip .2cm 

Supposons d'abord que $\Lambda _{\mu}(\mu _{[x],k})$ a un plus grand \'el\'ement $\gamma$ et choisissons un polyn\^ome $\phi$ dans $\Phi _{\mu}(\mu _{[x],k})$, tel que $\mu (\phi) = \gamma$. 
Alors si nous notons comme pr\'ec\'edement $\mu _{\phi ,\gamma}$ la valuation augment\'ee pour la valuation $\mu _k$ associ\'ee au polyn\^ome $\phi$ et \`a la valeur $\gamma$, tout polyn\^ome $f$ de $K[x]$ de degr\'e $deg (f) \leq d_{>k}$ v\'erifie $\mu _{\phi ,\gamma} (f) = \mu (f)$. 
Nous en d\'eduisons que la valuation $\mu _k$ appartient \`a une sous-famille discr\`ete $\mathcal{D}_{[x]}^{(j)}$ de la famille admise $\mathcal{A}_{[x]}$, et n'est pas \'egale \`a la derni\`ere valuation de cette sous-famille, l'\'el\'ement $k$ admet un successeur $k+1$ dans l'ensemble d'indices $I_{[x]}$, la valuation augment\'ee $\mu _{\phi ,\gamma}$ est la valuation $\mu _{[x],k+1}$ de la famille admise $\mathcal{A}_{[x]}$ et elle appartient encore \`a la sous-famille discr\`ete $\mathcal{D}_{[x]}^{(j)}$, nous notons $\phi _{k+1}$ le polyn\^ome-cl\'e $\phi$ et nous avons deg$\phi _{k+1} =d_{[x],k+1} = d_{>k}$.   
 
D'apr\`es le lemme \ref{le-bij} nous en d\'eduisons que $\Lambda _{\mu}(\mu _{[y],k}) $ a un plus grand \'el\'ement $\delta$, que le polyn\^ome $\psi = Inv_{x/y}(\phi )$ de $\Phi _{\mu}(\mu _{[y],k})$ v\'erifie $\mu (\psi ) = \delta$ et que tout polyn\^ome $g$ de $K[y]$ de degr\'e $deg (g) \leq d_{>k}$ v\'erifie $\mu _{\psi ,\delta} (g) = \mu (g)$.  
 De fa\c con analogue la valuation $\mu _k$ appartient \`a une sous-famille discr\`ete $\mathcal{D}_{[y]}^{(j)}$ de la famille admise $\mathcal{A}_{[y]}$, n'est pas \'egale \`a la derni\`ere valuation de cette sous-famille, l'\'el\'ement $k$ admet un successeur $k+1$ dans l'ensemble d'indices $I_{[y]}$, la valuation augment\'ee $\mu _{\psi ,\delta}$ est la valuation $\mu _{[y],k+1}$, et nous notons $\psi _{k+1}$ le polyn\^ome-cl\'e $\psi$ de degr\'e $d_{[y],k+1} = d_{>k}$.  

Par cons\'equent les ensembles $\{ j \in I_{[x]}\ | \ j \leq k+1\}$ et $\{ j \in I_{[y]}\ | \ j \leq k+1\}$ sont identiques, et d'apr\`es le lemme pr\'ec\'edent les familles ${\mathcal A}_{[x]} = \bigl ( \mu _{[x],k} \bigr )_{i \in I_{[x]}}$ et ${\mathcal A}_{[y]}= \bigl ( \mu _{[y],k} \bigr )_{i \in I_{[y]}}$ sont les m\^emes jusqu'\`a l'indice $k+1$, et nous posons $\mu _{k+1} = \mu _{[x],k+1} = \mu _{[y],k+1}$. 

\vskip .2cm 

Supposons maintenant que $\Lambda _{\mu}(\mu _{[x],k})$ n'a pas de plus grand \'el\'ement et que nous avons l'in\'egalit\'e stricte $d_{>k} > d_k$, alors nous sommes dans le cas o\`u la valuation $\mu _k$ est l'avant-derni\`ere valuation d'une sous-famille admissible discr\`ete $\mathcal{D}_{[x]}^{(j)}$ de la famille admise $\mathcal{A}_{[x]}$.  
Nous choisissons une valeur $\gamma _{k+1}$ dans $\Lambda _{\mu}(\mu _{[x],k})$ et un polyn\^ome $\phi _{k+1}$ dans $\Phi _{\mu}(\mu _{[x],k})$ tel que $\mu (\phi _{k+1}) = \gamma _{k+1}$ et nous d\'efinissons comme pr\'ec\'edemment la valuation $\mu _{[x],k+1}$ comme la valuation augment\'ee $\mu _{\phi _{k+1}, \gamma _{k+1}}$. 
Cette valuation est alors la derni\`ere valuation de la sous-famille $\mathcal{D}_{[x]}^{(j)}$,  et nous avons deg$\phi _{k+1} =d_{[x],k+1} = d_{>k}$.   

Nous d\'eduisons du lemme \ref{le-bij} que $\Lambda _{\mu}(\mu _{[y],k})$ n'a pas non plus de plus grand \'el\'ement, alors la valuation $\mu _k$ est aussi l'avant-derni\`ere valuation d'une sous-famille admissible discr\`ete $\mathcal{D}_{[y]}^{(j)}$ de la famille admise $\mathcal{A}_{[y]}$, et de la m\^eme mani\`ere que pr\'ec\'edemment, si nous choisissons le polyn\^ome $\psi _{k+1}$ dans $\Phi _{\mu}(\mu _{[y],k})$ \'egal \`a $Inv_{x/y}(\phi _{k+1})$ et d\'efinissons la valuation $\mu _{[y],k+1}$ comme la valuation augment\'ee $\mu _{\psi _{k+1}, \delta_{k+1}}$, avec $\delta _{k+1}= \mu (\psi _{k+1})$, les valuations $\mu _{[x],k+1}$ et $\mu _{[y],k+1}$ sont \'egales. 
Par cons\'equent les ensembles $\{ j \in I_{[x]}\ | \ j \leq k+1\}$ et $\{ j \in I_{[y]}\ | \ j \leq k+1\}$ sont identiques, et d'apr\`es le lemme pr\'ec\'edent les familles ${\mathcal A}_{[x]} = \bigl ( \mu _{[x],k} \bigr )_{i \in I_{[x]}}$ et ${\mathcal A}_{[y]}= \bigl ( \mu _{[y],k} \bigr )_{i \in I_{[y]}}$ sont les m\^emes jusqu'\`a l'indice $k+1$, et nous posons $\mu _{k+1} = \mu _{[x],k+1} = \mu _{[y],k+1}$. 

\vskip .2cm 

Supposons enfin que $\Lambda _{\mu}(\mu _{[x],k})$ n'a pas de plus grand \'el\'ement et que nous avons l'\'egalit\'e stricte $d_{>k} = d_k$, alors nous pouvons supposer que la valuation $\mu _k$ est la derni\`ere valuation d'une sous-famille admissible discr\`ete $\mathcal{D}_{[x]}^{(j)}$ de la famille $\mathcal{A}_{[x]}$. 

Il existe une famille $(\phi _{\alpha})_{\alpha \in A}$ de polyn\^omes appartenant \`a $\Phi _{\mu}(\mu _{[x],k})$ index\'ee par un ensemble totalement ordonn\'e $A$, sans plus grand \'el\'ement, telle que la famille des valeurs $(\gamma _{\alpha})_{\alpha \in A}$, o\`u $\gamma _{\alpha} = \mu (\phi _{\alpha})$ est cofinale dans  $\Lambda _{\mu}(\mu _{[x],k})$.  
La famille de valuations augment\'ees associ\'ee $(\mu _{\phi _{\alpha} , \gamma _{\alpha}})_{\alpha \in A}$ est une famille admissible continue de valuations de $K[x]$, et en particulier pour tout $\alpha < \beta$ dans $A$ le polyn\^ome $\phi _{\beta}$ appartient \`a l'ensemble $\Phi _{\mu}(\mu _{\phi _{\alpha} , \gamma _{\alpha}})$. 
Nous posons $A_{[x]}^{(j)} =A$, et pour tout $\alpha$ nous posons $\mu _{[x],\alpha} ^{(j)}= \mu _{\phi _{\alpha} , \gamma _{\alpha} }$, alors la famille $\mathcal{C}_{[x]}^{(j)} = {\bigl (\mu _{[x],\alpha}^{(j)}\bigr )}_{\alpha \in A_{[x]}^{(j)}}$ est la partie continue de la sous-famille admissible simple $\mathcal{S}_{[x]}^{(j)} $ de la famille admise $\mathcal{A}_{[x]}$. 

De ce qui pr\'ec\`ede nous en d\'eduisons que la valuation $\mu _k$ est aussi la derni\`ere valuation d'une sous-famille admissible discr\`ete $\mathcal{D}_{[y]}^{(j)}$ de la famille $\mathcal{A}_{[y]}$, et que la famille $\mathcal{C}_{[y]}^{(j)} = {\bigl (\mu _{[y],\alpha}^{(j)}\bigr )}_{\alpha \in A_{[y]}^{(j)}}$ est la partie continue de la sous-famille admissible simple $\mathcal{S}_{[y]}^{(j)} $ de la famille admise $\mathcal{A}_{[y]}$. 
Si nous posons encore $A_{[y]}^{(j)} =A$ et $\mu _{[y],\alpha} ^{(j)}= \mu _{\psi _{\alpha} , \delta _{\alpha} }$, avec $\psi _{\alpha} = Inv_{x/y}(\phi _{\alpha})$ et $\delta _{\alpha} = \mu (\psi _{\alpha})$, les valuations $\mu _{[x],\alpha} ^{(j)}$ et $\mu _{[y],\alpha} ^{(j)}$ sont \'egales. 

Par cons\'equent les $j$ premi\`eres sous-familles admissibles simples des familles admises $\mathcal{A}_{[x]}$ et $\mathcal{A}_{[y]}$ co\"\i ncident et il reste \`a montrer que si la valuation $\mu$ est diff\'erente de la valuation limite de la famille continue $\mathcal{C}^{(j)} = {\bigl (\mu _{\alpha}^{(j)}\bigr )}_{\alpha \in A^{(j)}}$, les premi\`eres valuations des sous-familles admissibles respectivement $\mathcal{S}_{[x]}^{(j+1)}$ et $\mathcal{S}_{[y]}^{(j+1)}$ sont \'egales. 
 
Comme ces valuations $\mu _{[x],1} ^{(j+1)}$ et $\mu _{[y],1} ^{(j+1)}$ sont obtenues comme valuations augment\'ees limites de la famille continue $\mathcal{C}^{(j)}$ nous devons consid\'erer les ensembles: 

$$ \tilde\Phi _{[x]}\bigl ( \mathcal{C} ^{(j)} \bigr ) \ = \ \{ f \in K[x] \ | \ \mu _{\alpha} ^{(j)} (f) < \mu _{\beta}^{(j)} (f) , \forall \alpha < \beta \in A^{(j)} \}$$ 
et 
$$\tilde\Phi _{[y]}\bigl ( \mathcal{C} ^{(j)} \bigr )  \ = \ \{ g \in K[y] \ | \ \mu _{\alpha} ^{(j)} (g) < \mu _{\beta}^{(j)} (g) , \forall \alpha < \beta \in A^{(j)} \}$$ 
 
Par hypoth\`ese ces ensembles sont non vides, nous appelons respectivement $d _{[x], \mathcal{C} ^{(j)}}$ et $d _{[y], \mathcal{C} ^{(j)}}$ le degr\'e minimal des polyn\^omes appartenant \`a $\tilde\Phi _{[x]} \bigl ( \mathcal{C} ^{(j)} \bigr )$ et \`a $\tilde\Phi _{[y]} \bigl ( \mathcal{C} ^{(j)} \bigr )$, et nous d\'efinissons les ensembles  
$$\Phi _{[x]}\bigl ( \mathcal{C} ^{(j)} \bigr ) \ = \ \{ \phi \in \tilde\Phi _{[x]}\bigl ( \mathcal{C} ^{(j)} \bigr ),\  \hbox{deg}\phi =  d _{[x], \mathcal{C} ^{(j)}} \ \hbox{et $\phi$ unitaire} \  \} $$  
et 
$$\Phi _{[y]}\bigl ( \mathcal{C} ^{(j)} \bigr ) \ = \ \{ \psi \in \tilde\Phi _{[y]}\bigl ( \mathcal{C} ^{(j)} \bigr ),\  \hbox{deg}\psi =  d _{[y], \mathcal{C} ^{(j)}} \ \hbox{et $\psi$ unitaire} \  \} \ .$$   
\vskip .2cm 

\begin{lemma} \label{le-deg}
Nous avons l'\'egalit\'e $d _{[x], \mathcal{C} ^{(j)}} = d _{[y], \mathcal{C} ^{(j)}}$ et l'application $Inv_{x/y}$ induit une bijection de $\Phi _{[x]}\bigl ( \mathcal{C} ^{(j)} \bigr )$ sur $\Phi _{[y]}\bigl ( \mathcal{C} ^{(j)} \bigr )$. 
\end{lemma}

\begin{preuvele} 
Comme dans la preuve du lemme \ref{le-bij} nous v\'erifions que si un polyn\^ome $f$ v\'erifie $\mu _{\alpha} ^{(j)}(f) < \mu _{\beta} ^{(j)}(f)$ pour tout $\alpha < \beta$ dans $A ^{(j)}$, alors le polyn\^ome $g = Inv_{x/y}(f)$ est de m\^eme degr\'e que $f$ et v\'erifie les m\^emes in\'egalit\'es. 

Par cons\'equent l'application $Inv_{x/y}$ induit une bijection de $\tilde\Phi _{[x]}\bigl ( \mathcal{C} ^{(j)} \bigr ) \cap IU(K[x])$ sur $\tilde\Phi _{[y]}\bigl ( \mathcal{C} ^{(j)} \bigr ) \cap IU(K[y])$, d'o\`u l'\'egalit\'e des degr\'es $d _{[x], \mathcal{C} ^{(j)}} = d _{[y], \mathcal{C} ^{(j)}}$ et la bijection cherch\'ee. 

\hfill \tf    
\end{preuvele} 
 
\begin{lemma} \label{le-egal-limite}
Soient $\phi$ et $\psi = Inv_{x/y}(\phi )$ deux polyn\^omes appartenant respectivement aux ensembles $\Phi _{[x]}\bigl ( \mathcal{C} ^{(j)} \bigr )$ et $\Phi _{[y]}\bigl ( \mathcal{C} ^{(j)} \bigr )$, et soient $\gamma = \mu (\phi )$ et $\delta = \mu (\psi)$. 
Alors les valuations augment\'ees limites $\mu ' _{[x]} = \bigl [ \mathcal{C} ^{(j)} \ ; \ \mu ' _{[x]} (\phi ) = \gamma \bigr ]$ et $\mu ' _{[y]}= \bigl [ \mathcal{C} ^{(j)} \ ; \ \mu ' _{[y]}(\psi ) = \delta \bigr ]$ sont \'egales. 
\end{lemma} 

\begin{preuvele} 
La preuve est analogue \`a la preuve du lemme \ref{le-valuations-egales}, nous utilisons la proposition \ref{prop-definition-val-augm-lim} \`a la place du th\'eor\`eme 5.2. \cite{McL 1}. 

\hfill \tf    
\end{preuvele} 
 
Pour conna\^\i tre les premi\`eres valuations $\mu _{[x],1} ^{(j+1)}$ et $\mu _{[y],1} ^{(j+1)}$ des sous-familles admissibles $\mathcal{S}_{[x]}^{(j+1)}$ et $\mathcal{S}_{[y]}^{(j+1)}$, il faut consid\'erer les ensembles $\Lambda_{[x]}\bigl ( \mathcal{C} ^{(j)} \bigr )$ et $\Lambda _{[y]}\bigl ( \mathcal{C} ^{(j)} \bigr )$ d\'efinis respectivement par  
$$ \Lambda _{[x]}\bigl ( \mathcal{C} ^{(j)} \bigr ) = \{ \mu (\phi ) \ | \ \phi \in \Phi _{[x]}\bigl ( \mathcal{C} ^{(j)} \bigr ) \} \ \hbox{et} \  \Lambda _{[y]}\bigl ( \mathcal{C} ^{(j)} \bigr ) = \{ \mu (\psi ) \ | \ \psi \in \Phi _{[y]}\bigl ( \mathcal{C} ^{(j)} \bigr ) \} \ .$$ 

Nous d\'eduisons du lemme \ref{le-deg} que l'ensemble $\Lambda_{[x]}\bigl ( \mathcal{C} ^{(j)} \bigr )$ a un plus grand \'el\'ement $\gamma$ si et seulement si $\Lambda_{[y]}\bigl ( \mathcal{C} ^{(j)} \bigr )$ a un plus grand \'el\'ement $\delta$, et que si $\phi$ est un polyn\^ome de $\Phi _{[x]}\bigl ( \mathcal{C} ^{(j)} \bigr )$ v\'erifiant $\mu (\phi ) = \gamma$ alors le polyn\^ome $\psi$ de $\Phi _{[y]}\bigl ( \mathcal{C} ^{(j)} \bigr )$ d\'efini par $\psi = Inv _{x/y}(\phi )$  v\'erifie $\mu (\psi ) = \delta$. 

Et de m\^eme si $\Lambda_{[x]}\bigl ( \mathcal{C} ^{(j)} \bigr )$ n'a pas de plus grand \'el\'ement et si $\bigl ( \gamma _{\theta} \bigl ) _{\theta \in T}$ est une famille cofinale dans $\Lambda_{[x]}\bigl ( \mathcal{C} ^{(j)} \bigr )$ avec $\gamma _{\theta} = \mu ( \phi _{\theta} )$, la famille $\bigl ( \delta _{\theta} \bigl ) _{\theta \in T}$ d\'efinie par $\delta _{\theta} = \mu ( \psi _{\theta} )$, avec $\psi _{\theta} = Inv _{x/y} ( \phi _{\theta} )$ est aussi une famille cofinale dans $\Lambda_{[y]}\bigl ( \mathcal{C} ^{(j)} \bigr )$.  
 
Dans le cas o\`u $\Lambda_{[x]}\bigl ( \mathcal{C} ^{(j)} \bigr )$  et $\Lambda_{[y]}\bigl ( \mathcal{C} ^{(j)} \bigr )$ ont un plus grand \'el\'ement, respectivement $\gamma = \mu (\phi )$ et $\delta = \mu (\psi )$, nous posons $\phi _{[x],1} ^{(j+1)} = \phi$, $\gamma _{[x],1} ^{(j+1)} = \gamma$ et  $\psi _{[y],1} ^{(j+1)} = \psi$, $\delta _{[y],1} ^{(j+1)} = \delta$. Alors les valuation $\mu _{[x],1} ^{(j+1)}$ et $\mu _{[x],1} ^{(j+1)}$ sont les valuations augment\'ees limites $\mu _{[x],1} ^{(j+1)} = \bigl [ \mathcal{C} ^{(j)} \ ; \ \mu _{[x],1} ^{(j+1)} ( \phi _{[x],1} ^{(j+1)}) = \gamma _{[x],1} ^{(j+1)} \bigr ]$ et $\mu _{[y],1} ^{(j+1)} = \bigl [ \mathcal{C} ^{(j)} \ ; \ \mu _{[y],1} ^{(j+1)} ( \psi _{[y],1} ^{(j+1)}) = \delta _{[y],1} ^{(j+1)} \bigr ]$,  et nous d\'eduisons du lemme \ref{le-egal-limite} que ces valuations sont \'egales. 

Dans le cas o\`u les ensembles $\Lambda_{[x]}\bigl ( \mathcal{C} ^{(j)} \bigr )$  et $\Lambda_{[y]}\bigl ( \mathcal{C} ^{(j)} \bigr )$ n'ont pas de plus grand \'el\'ement, nous pouvons choisir un polyn\^ome quelconque $\phi _{[x],1} ^{(j+1)}$ dans $\Phi _{[x]} \bigl ( \mathcal{C} \bigr )$ et son image $\psi _{[y],1} ^{(j+1)} = Inv_{x/y} (\phi _{[x],1} ^{(j+1)})$ dans $\Phi _{[y]} \bigl ( \mathcal{C} \bigr )$ pour d\'efinir les valuations $\mu _{[x],1} ^{(j+1)}$ et $\mu _{[x],1} ^{(j+1)}$ comme les valuations augment\'ees limites $\mu _{[x],1} ^{(j+1)} = \bigl [ \mathcal{C} ^{(j)} \ ; \ \mu _{[x],1} ^{(j+1)} ( \phi _{[x],1} ^{(j+1)}) = \gamma _{[x],1} ^{(j+1)} \bigr ]$ et $\mu _{[y],1} ^{(j+1)} = \bigl [ \mathcal{C} ^{(j)} \ ; \ \mu _{[y],1} ^{(j+1)} ( \psi _{[y],1} ^{(j+1)}) = \delta _{[y],1} ^{(j+1)} \bigr ]$,  et nous d\'eduisons encore du lemme \ref{le-egal-limite} que ces valuations sont \'egales. 

\vskip .2cm 

Nous avons ainsi montr\'e que si $x$ et $y$ sont deux g\'en\'erateurs de l'extension $L$ sur $K$ nous pouvons choisir des familles admises $\mathcal{A}_{[x]}$ et $\mathcal{A}_{[y]}$ qui sont \'egales en tant que familles de valuations du corps $L$. 

\hfill \tf    
\end{preuve}   
\vskip .2cm 

Ainsi \`a toute valuation $\mu$ du corps $L$ extension transcendante pure de $K$, prolongeant la valuation $\nu$, nous pouvons associer une famille $\mathcal{A} (\mu ) = \left ( \mu _i \right ) _{i \in I}$ de valuations de $L$, que nous appelons encore \emph{la famille admise} associ\'ee \`a la valuation $\mu$. 

Comme nous l'avons remarqu\'e la notion de famille admissible ou admise de valuations de l'anneau des polyn\^omes $K[x]$ ne peut pas s'\'etendre simplement \`a une famille de valuations du corps $L=K(x)$. 
Nous pouvons donc nous demander s'il est possible de donner une d\'efinition de famille admissible ou de famille admise de valuations de $L$ qui ne fasse pas intervenir le choix d'un g\'en\'erateur. 

\vskip .2cm 

Dans le cas d'une extension alg\'ebrique $L=K(\theta )$ et d'une pseudo-valuation $\mu$ de $K[x]$ de socle l'id\'eal engendr\'e par le polyn\^ome minimal de $\theta$ sur $K$, il ne semble pas possible de trouver une famille de valuations de $K[x]$ in d\'ependante du g\'en\'erateur $\theta$. 

Mais par analogie avec l'\'etude des singularit\'es des courbes planes, il semble naturel de se demander si au moins pour des g\'en\'erateurs suffisamment g\'en\'raux, dans un sens \`a pr\'eciser, de l'extension $L$, les familles admises associ\'ees partagent certaines propri\'et\'es. 
Par exemple nous pouvons nous demander si le nombre de sous-familles admissibles simples dans chaque famille admise est ind\'ependant du g\'en\'erateur, ou si le nombre de valuations dans chaque sous-famille admissible discr\`ete est le m\^eme.

		  

\vskip 1cm


\begin{thebibliography}{xxxxx}
%

\bibitem [McL 1] {McL 1} S. MacLane: A construction for absolute values in polynomial rings. Trans. Amer. Math. Soc. {\bf 40} (1936), 363-395.

\bibitem [McL 2] {McL 2} S. MacLane: A construction for prime ideals as absolute values of an algebraic field. Duke Math. J. {\bf 2} (1936), 492-510.

\bibitem [Va1] {Va 1} M. Vaqui\'e: Extension d'une valuation. Trans. Amer. Math. Soc. {\bf 359} (2007), 3439-3481.

\bibitem [Va2] {Va 2} M. Vaqui\'e: Famille admise associ\'ee \`a une valuation de $K[x]$, Sem. et Congr. {\bf 10} (2005), 391-428. 

%
                    \end{thebibliography}
\end{document}